\documentclass[aos,preprint,tikz,border=2pt,png]{imsart}
\pdfoutput = 0

\arxiv{arXiv:0000.0000}
\usepackage[toc]{appendix}
\usepackage{tkz-euclide}
\usetkzobj{all}

\usepackage[hyphens]{url}
\RequirePackage[OT1]{fontenc}
\RequirePackage{amsthm,amsmath, amsfonts, amssymb, graphics, graphicx, subfigure, color, multirow}

\RequirePackage{natbib}
\RequirePackage[colorlinks,citecolor=blue,urlcolor=blue,breaklinks]{hyperref}

\usepackage{pstricks-add}
\usepackage[draft]{todonotes}
\usepackage{algorithm}
\usepackage{algpseudocode}

\newif\ifdraft
\draftfalse 

\usepackage{graphicx}
\usepackage[singlelinecheck=false]{caption}

\usepackage{algpseudocode}
\usepackage{comment}

\newtheorem{lem}{Lemma}
\newtheorem{cor}{Corollary}

\newtheorem{thm}{Theorem}

\newtheorem{definition}{Definition}
\newtheorem{rmk}{Remark}
\usepackage{titlesec}

\newcommand{\vxbarbar}{ \overline{\overline{\vx}} }

\newcommand{\iid}{\stackrel{\mathrm{iid}}{\sim}}

\newcommand{\vM}{\boldsymbol{M}}

\newcommand{\vm}{\boldsymbol{m}}

\newcommand{\bbE}{\mathbb{E}}
\newcommand{\bbP}{\mathbb{P}}

\newcommand{\vA}{\boldsymbol{A}}

\newcommand{\vB}{\boldsymbol{B}}

\newcommand{\vE}{\boldsymbol{E}}

\newcommand{\vI}{\mathbf{I}}

\newcommand{\vP}{\boldsymbol{P}}
\newcommand{\vp}{\boldsymbol{p}}

\newcommand{\vu}{\boldsymbol{u}}
\newcommand{\vV}{\boldsymbol{V}}
\newcommand{\vv}{\boldsymbol{v}}
\newcommand{\vx}{\boldsymbol{x}}
\newcommand{\vX}{\boldsymbol{X}}

\newcommand{\vw}{\boldsymbol{w}}
\newcommand{\vy}{\boldsymbol{y}}

\newcommand{\vz}{\boldsymbol{z}}

\newcommand{\vT}{\boldsymbol{T}}

\newcommand{\vzero}{\boldsymbol{0}}

\newcommand{\bbeta}{\boldsymbol{\beta}}

\newcommand{\bgamma}{\boldsymbol{\gamma}}

\newcommand{\btheta}{\boldsymbol{\theta}}

\newcommand{\bLambda}{\boldsymbol{\Lambda}}

\newcommand{\bSigma}{\boldsymbol{\Sigma}}

\newcommand{\bepsilon}{\boldsymbol{\epsilon}}

\newcommand{\boldeta}{\boldsymbol{\eta}}
\newcommand{\bxi}{\boldsymbol{\xi}}

\newcommand\independent{\protect\mathpalette{\protect\independenT}{\perp}}
\def\independenT#1#2{\mathrel{\rlap{$#1#2$}\mkern2mu{#1#2}}}

\newcommand{\epf}{\hfill $\Box$}

\newcommand*{\Scale}[2][4]{\scalebox{#1}{$#2$}}%

\begin{document}

\begin{frontmatter}

\title{On consistency and sparsity for sliced inverse regression in high dimensions}
\runtitle{DT-SIR }
\thankstext{t1}{Lin's research is supported by the Center of Mathematical Sciences and Applications at Harvard University. Zhao's research is supported by the NSF Grant DMS-1208735. Liu's research is supported by the NSF Grant DMS-1120368 and NIH  Grant R01 GM113242-01}

\begin{aug}
 \author{\fnms{Qian} \snm{Lin}\thanksref{m1}\ead[label=e1]{qianlin@cmsa.fas.harvard.edu}\thanksref{t1}}
 \author{\fnms{Zhigen} \snm{Zhao}\thanksref{m2}\ead[label=e2]{zhaozhg@temple.edu}\thanksref{t1}}
\and    
\author{\fnms{Jun S.} \snm{Liu}\thanksref{m1}\ead[label=e3]{jliu@stat.harvard.edu}\thanksref{t1}}

\runauthor{Q. Lin, Z. Zhao and J. S. Liu}

\affiliation{Harvard University\thanksmark{m1} }\affiliation{Temple University\thanksmark{m2} }
  
  \address{Qian Lin\\
    Center of Mathematical Sciences\\
     and Applications \\
    Harvard University\\
    One Oxford Street \\
    Cambridge, MA 02138  \\
    USA\\
    \printead{e1}}    
  \address{Zhigen Zhao\\
    Department of Statistics \\
    Temple University\\
    346 Speakman Hall \\
    1810 N. 13th Street \\
    Philadelphia, PENNSYLVANIA, 19122 \\
    USA\\
    \printead{e2}}   
  \address{Jun S.  Liu\\
    Department of Statistics\\
    Harvard University\\
    1 Oxford Street \\
    Cambridge, MA 02138  \\
    USA\\
    \printead{e3}}  
\end{aug}

\begin{abstract}

We provide here a framework to analyze the phase transition phenomenon of slice inverse regression (SIR), a supervised dimension reduction technique introduced by \cite{Li:1991}.  Under mild conditions, the asymptotic ratio $\rho= \lim p/n$ is the phase transition parameter and the SIR estimator is consistent if and only if  $\rho= 0$. When dimension $p$ is greater than $n$, we propose a diagonal thresholding screening SIR (DT-SIR) algorithm. This 
method provides us with an estimate of the eigen-space of the covariance matrix of the conditional expectation $var(\bbE[\vx|y])$. The desired dimension reduction space is then obtained by multiplying the inverse of the covariance matrix on the eigen-space. Under certain sparsity assumptions on both the covariance matrix of predictors and the loadings of the directions, we prove the consistency of DT-SIR  in estimating the dimension reduction space in high dimensional data analysis. Extensive numerical experiments demonstrate superior performances of the proposed method in comparison to its competitors.

\end{abstract}

\begin{keyword}[class=MSC]
\kwd[Primary ]{62J02}
\kwd[; secondary ]{62H25}
\end{keyword}

 \begin{keyword}
 \kwd{dimension reduction}
 \kwd{random matrix theory}
 \kwd{sliced inverse regression}
 \end{keyword}

\end{frontmatter}

\section{Introduction.}\label{sec:intro}

For a continuous multivariate random variable $(y,\vx)$ where $\vx \in \mathbb{R}^{p}$ and $y$ $\in \mathbb{R}$, a subspace $\mathcal{S}'\subset \mathbb{R}^p$ is called the effective dimension reduction (EDR) space if $y\independent  \vx | P_{\mathcal{S}'}(\vx)$ where $\independent$ stands for independence.  Under mild conditions (\cite{Cook:1996}), the intersection of all the EDR spaces is again an EDR space, which is denoted as $\mathcal{S}$ and called the central space.
Many algorithms were proposed to find such subspace $\mathcal{S}$ under the assumption  $d=dim\mathcal{S}\ll p$. 
This line of research is commonly known as sufficient dimension reduction. The Sliced Inverse Regression (SIR, \cite{Li:1991}) is the first, yet the most widely used method in sufficient dimension reduction, due to its simplicity, computational efficiency and generality. The asymptotic properties of SIR are of particular interest in the last two decades. 
The consistency of SIR has been proved  for fixed $p$ in \cite{Li:1991}, \cite{hsing1992asymptotic}, \cite{Zhu:Ng:1995} and \cite{Zhu:Fang:1996}. 
Later, \cite{Zhu:Miao:Peng:2006} have obtained the consistency if  $p=o(\sqrt{n})$. 
A similar restriction also appears in two recent work (see  \cite{zhong2012correlation} and \cite{jiang2013sliced}).
When $p>n$, a common strategy pursued by many recent researchers is to make sparsity assumptions that only a few predictors play a role in explaining and predicting $y$ and apply various regularization methods. 
For instance, \cite{li2006sparse, li2007sparse} and \cite{yu2013dimension} applied LASSO (\cite{tibshirani1996regression}), Dantzig selector (\cite{candes2007dantzig}) and elastic net (\cite{zou2005regularization}) respectively to solve the generalized eigenvalue problems raised by a variety of SDR algorithms.

However, a piece of jigsaw is missing in the understanding of SIR. If the dimension $p$ diverges as $n$ increases, when will the SIR break down? A similar question has been asked for a variety of SDR estimates in \cite{cook2012estimating}.  
 In this paper, we prove that, under certain technical assumptions, the SIR estimator is consistent if and only if $\rho=\lim\frac{p}{n}=0$. 
Such a result on inconsistency provides theoretical justifications for imposing certain structural assumption, such as sparsity, in high dimensional settings. 
This behavior of SIR in high dimension, which will be called the phase transition phenomenon, is similar to that of the principal component analysis (PCA), an unsupervised counterpart of SIR. 
This extension is, however, by no means trivial. 
 After  all the samples $(y_{i},\vx_{i})$ are sliced into H bins according to the order statistics of $y_{i}$ , the sliced samples are neither independent nor identically distributed. This difference  increases the difficulty significantly.
In this paper, we provide a new framework to study the phase transition behaviour of SIR. 
The technical tools developed here can potentially be extended to study the phase transition behaviour of other SDR estimators. 

The second part of the article aims at extending the original SIR to the scenario with ultra-high dimension $(p=o(\exp(n^{\xi})))$.
Based on  equation \eqref{eta:beta} in Section \ref{sec:notation}, the central space can be estimated by the column space of $\Scale[0.8]{\widehat{\bSigma}_{\vx}^{-1}}col(\Scale[0.8]{\widehat{\vV}_{H}})$, where $\Scale[0.8]{\widehat{\bSigma}_{\vx}^{-1}}$ is any consistent estimate of the precision matrix  $\Scale[0.8]{\bSigma_{\vx}^{-1}}$ and  $col(\Scale[0.8]{\widehat{\vV}_{H}})$ is the estimate of the space $col(var(\bbE[\vx|y]))$.  To estimate the column space of  $var\left(\bbE[\vx|y]\right)$,  we propose a diagonal screening procedure based on new univariate statistics $var_{H}(\vx(k))$, which are the diagonal elements of $var\left(\bbE[\vx|y]\right)$, motivated by recent work in sparse PCA (\cite{johnstone2004sparse}).
 After ranking the predictors according to the magnitude of $var_{H}(\vx(k))$ decreasingly, we choose the set $\mathcal{I}$ consisting of the first $R$ predictors as active predictors.
 The SIR procedure is subsequently applied to these selected predictors to estimate  the  $d$-dimensional column space of $var(\bbE[\vx^{\mathcal{I}}|y])$ by $col(\Scale[0.8]{\widehat{\vV}^{\mathcal{I}}_{H}})$   where $\Scale[0.8]{\widehat{\vV}^{\mathcal{I}}_{H}}$ is the matrix formed by the top $d$ eigenvectors 
 of $\Scale[0.8]{\widehat{\bLambda}^{\mathcal{I}}_{H}}$.
We embed $\Scale[0.8]{\widehat{\vV}^{\mathcal{I}}_{H}}$ into $\mathbb{R}^{p\times d}$ by filling in 0's for entries outside the chosen row set $\mathcal{I}$, and denote this new matrix by $e(\Scale[0.8]{\widehat{\vV}^{\mathcal{I}}_{H}})$. 
The estimate of the central space is defined to be $col(\Scale[0.8]{\widehat{\bSigma}_{\vx}^{-1}}e(\Scale[0.8]{\widehat{\vV}^{\mathcal{I}}_{H})})$. We  name this two-stage algorithm  as {\bf D}iagonal {\bf T}hresholding  {\bf SIR} (DT-SIR), and
 prove that DT-SIR is consistent in estimating the central space under certain regularity conditions. Extensive simulation studies show that DT-SIR performs better than its competitors and is computationally efficient.

The rest of the paper is organized as follows. In Section \ref{sec:notation}, we  briefly describe the SIR procedure and introduce the notations. 
In Section \ref{sec:sir}, after a brief review of existing asymptotic results of SIR 
procedure, we state Theorems \ref{thm:space:consistency} and \ref{thm:main1.1} to discuss the phase transition phenomenon of SIR. 
In Section \ref{sec:high_dimensional}, we propose the DT-SIR method and show that DT-SIR is consistent in high dimensional data analysis. In Section \ref{sec:simulation}, we provide simulation studies to compare DT-SIR with its competitors. Concluding remarks and discussions are put in Section \ref{sec:conclusion}. All the proofs are presented in appendices.

\section{Preliminaries and notations.}\label{sec:notation}
\subsection{Sliced inverse regression}
Consider the multiple index model
\begin{equation}\label{model:li:1991}
y=f(\bbeta_{1}^{\tau}\vx,\cdots,\bbeta_{d}^{\tau}\vx,\epsilon)
\end{equation}
where $\vx\in \mathbb{R}^{p}$, $\epsilon$ is the noise and $f$ is an unknown link function. 
Without loss of generality,  we assume that $\mathbb{E}[\vx]=0 \in \mathbb{R}^{p}$.
Although the $p\times d$ matrix $\vV=(\bbeta_{1},\cdots,\bbeta_{d})$ is not identifiable, the space spanned by the $\bbeta$'s, which is called the column space of $\vV$ and denoted by $col(\vV)$,  might be identified.
\cite{Li:1991} proposed the {\it Sliced Inverse Regression} (SIR) procedure to estimate the central space $col(\vV)$ without knowing $f(\cdot)$, which can be briefly summarized as follows: 
Given $n$ i.i.d. samples $(y_{i},\vx_{i})$, $i=1,\cdots,n$, SIR first divides them into $H$ equal-sized slices according to the order statistics $y_{(i)}$.\footnote{To ease notations and arguments, we assume that $n=cH$ and $H=o\left(\log(n)\wedge \log(p)\right)$ throughout the article.}
We re-express the data as $y_{h,j}$ and $\vx_{h,j}$, where $(h,j)$ is the double subscript in which $h$ refers to the slice number and $j$ refers to the order number of a sample in the $h$-th slice, i.e.,
\begin{equation*}
y_{h,j}=y_{(c(h-1)+j)}, \mbox{\quad \quad } \vx_{h,j}=\vx_{(c(h-1)+j)}.
\end{equation*}
Here $\vx_{(k)}$ is the concomitant of $y_{(k)}$.
Let the sample mean in the $h$-th slice 
be  $
\overline{\vx}_{h,\cdot}$, and let the mean of all the samples be $\vxbarbar$. 
 Then, $\bLambda_p\triangleq var(\bbE[\vx|y])$ can be estimated by:
\begin{equation}\label{eqn:lambda}
\widehat{\bLambda}_{H}=\frac{1}{H}\sum_{h=1}^{H}\bar{\vx}_{h,\cdot}\bar{\vx}_{h,\cdot}^{\tau}. 
\end{equation}
Based on the observation that
\begin{equation}\label{eta:beta}
col(\bLambda)=\bSigma_{x}col(\vV),
\end{equation}
the SIR then estimates the central space $col(\vV)$ by $\widehat{\Sigma}^{-1}_{\vx}col(\widehat{\vV}_{H})$  where $\Scale[0.8]{\widehat{\vV}_{H}}$ is the matrix formed by the top $d$ eigenvectors 
 of $\Scale[0.8]{\widehat{\bLambda}_{H}}$..  Throughout the article, we assume that $d$ is fixed and  the $d$-th largest eigenvalue $\lambda_{d}$ of $\bLambda_{p}$ is bounded away from 0 when $n, p\rightarrow \infty$.

In order  for SIR to result in a consistent estimate of the central space, \cite{Li:1991} imposed the 
the following two conditions:
\begin{itemize} \label{cond:linear}
\item ({\bf A1}).  {\bf Linearity condition:}
For any $\bxi \in \mathbb{R}^{p}$, 
$
\mathbb{E}[\bxi^{\tau}\vx|\bbeta_{1}^{\tau}\vx,\cdots,\bbeta_{d}^{\tau}\vx]$ is a linear combination of $\bbeta_{1}^{\tau}\vx,\cdots,\bbeta_{d}^{\tau}\vx.
$
\item ({\bf A2}). {\bf Coverage condition:} The dimension of the space spanned by the central curve equals the dimension of the central space, i.e., $d'=d$. 
\end{itemize}

\subsection{Further Notations.}
 Let $S_{h}$ be the  $h$-th interval $(y_{h-1,c},y_{h,c}]$ for $2\leq h \leq H-1$, $S_{1}=(-\infty, y_{1,c}]$ and $S_{H}=(y_{H-1,c},\infty)$. 
Note that these intervals depend on the order statistics $y_{(i)}$ and are thus random.
 For any $\omega$ in the product sample space, we define a random variable $\delta_h=\delta_h(\omega)=\int_{y\in S_h(\omega)}f(y)dy$ where $f(y)$ is the density function of $y$.  For $\mathcal{I} \subset \{1,\cdots,n\}, \mathcal{J} \subset\{1,\cdots,p\}$ and a $n\times p$ matrix $\vA$, $\vA^{\mathcal{I},\mathcal{J}}$ denotes the $|\mathcal{I}|\times |\mathcal{J}|$ sub-matrix formed by restricting the rows of $\vA$ to $\mathcal{I}$ and columns to $\mathcal{J}$. In articular, $\vA^{-,\mathcal{J}}$ denotes the sub-matrix formed by restricting the columns to $\mathcal{J}$;
 For a matrix $\vB=\vA^{\mathcal{I},\mathcal{J}} \in \mathbb{R}^{|\mathcal{I}|\times |\mathcal{J}|}$, we embed it into $\mathbb{R}^{p\times p}$ by putting 0 on entries outside $\mathcal{I}\times\mathcal{J}$ and denote the new matrix as $e(\vB)$.  Similar notations apply to vectors. 
For two positive numbers $a$ and $b$, we let $a\vee b \equiv \max\{a,b\}$ and let $a\wedge b\equiv \min\{a,b\}$.
 Let $\tau(x,t) = x \times 1(|x|>t)$ be the hard thresholding function.
 Throughout the article, $C$, $C_1$ and $C_2$ are used to denote generic absolute constants, though the actual value may vary from case to case.
For a vector $\vx$, we denote its $k$-th entry as $\vx(k)$.
 Let $\bbeta_1$ and $\bbeta_2$ be two vectors with the same dimension, the angle between these two vectors is denoted as $\angle(\bbeta_1,\bbeta_2)$.
 For two sequences $\{a_{n}\}$, $\{b_{n}\}$, we let $a_{n}\ll b_{n}$ stand for $a_{n}=O(b^{\epsilon}_{n})$ for some positive $\epsilon <1$ and let $a_{n}\succ b_{n}$ stand for $\lim\frac{b_{n}}{a_{n}}=0$.


\section{Consistency of SIR.}\label{sec:sir}

In order to control the behavior of SIR, we need to impose the following boundedness condition {\bf(A3)} on the predictors' covariance matrix in addition to the  tail condition (sub-Gaussian) on  their joint distribution. We also need a condition {\bf(A4)} for the central curve.

\begin{itemize}
\item {\bf (A3) Boundedness Condition: } $\vx$ is sub-Gaussian; and
there exist positive constants $C_{1}, C_{2}$ such that
\[
C_{1} \leq \lambda_{min}(\bSigma_{\vx})\leq \lambda_{max}(\bSigma_{\vx})\leq C_{2}
\]
where $\lambda_{min}(\bSigma_{\vx})$ and $\lambda_{max}(\bSigma_{\vx})$ are the minimal and maximal eigenvalues of  $\bSigma_{\vx}$ respectively.

 \item ({\bf A4})  The central curve $\vm(y)\triangleq\bbE[\vx|y]$ has finite fourth moment and is $\vartheta$-sliced stable (defined below) with respect to $y$ and $\vm(y)$.  
\end{itemize}

\begin{definition}\label{Definition:Sliced}
For two positive constants ${\bgamma}_{1}<1<\bgamma_{2}$, let $\mathcal{A}_{H}(\bgamma_{1},\bgamma_{2})$ be the collection of all the partition 
$-\infty=a_{0} < a_{1} < \cdots<a_{H-1} <a_{H}= \infty$
of $\mathbb{R}$ satisfying that
\[
\frac{\bgamma_{1}}{H} \leq P(a_{i}\leq y< a_{i+1} ) \leq \frac{\bgamma_{2}}{H}.
\]
The central curve $\vm(y)=\bbE[\vx|y]$  is called $\vartheta$-sliced stable with respect to $y$ for some $\vartheta>0$  if there exist positive constants $\bgamma_{i}, i=1,2,3$ such that for any $\bbeta$ in the central space for any partition in $\mathcal{A}_{H}(\bgamma_{1},\bgamma_{2})$, we have
\begin{equation}\label{def:sliced_stable}
\frac{1}{H}\Big|\sum_{h=0}^{H-1}var(\bbeta^{\tau}\vm(y)\mid a_{h}\leq y\leq a_{h+1})\Big|   \leq  \frac{\bgamma_{3}}{H^{\vartheta}}var(\bbeta^{\tau}\vm(y)).
\end{equation}
The central curve is sliced stable if it is $\vartheta$-sliced stable for some positive constant $\vartheta$.
\end{definition}


\begin{rmk} 
Note that we only need  \eqref{def:sliced_stable} to hold for all unit vectors in the central space by rescaling. 
By considering the orthogonal decomposition of $\bbeta$ in a general space with respect to the central space and its complement, it is easy to see that the sliced stability implies that \eqref{def:sliced_stable} holds true for all  vector $\bbeta$. In particular, we have the following two useful consequences of the slice-stability. 
\begin{itemize}
\item[i)]
By choosing $\bbeta^{\tau}=(0,\ldots,0,1,0,\ldots,0)$ with 1 at the k-th position, we have
\[
|\sum_{h=0}^{H}var(\vm(y,k)\mid a_{h}\leq y\leq a_{h+1})| \leq \bgamma_{3} H^{1-\vartheta} var(\vm(y,k)),
\]
where $\vm(y,k)$ is the k-th coordinate of the central curve $\vm(y)$.  \\

\item[ii)] Since equation \eqref{def:sliced_stable} holds for all unit vector $\bbeta$, we have
\[
\|\sum_{h=0}^{H}var(\vm(y)\mid a_{h}\leq y\leq a_{h+1})\|_{2} \leq \bgamma_{3} H^{1-\vartheta} \|var(\vm(y))\|_{2} .
\]
\end{itemize}
\end{rmk}

\begin{rmk}
  Suppose $\bbE[\vm(y)]=0$ and there are $n$ samples $\vm_{i}\triangleq\vm(y_{i})$. Let  $\vm_{h,i}$ and $\overline{\vm}_{h,\cdot}$ be defined similarly to $\vx_{h,i}$ and $\overline{\vx}_{h,\cdot}$, respectively.  On one hand,
 we have the classic consistent estimator $\frac{1}{n}\sum_{i}\vm_{i}\vm_{i}^{\tau}$ of $var(\vm(y))$.  
On the other hand,
if we expect that the slice-based estimate $\frac{1}{H}\sum_{h}\overline{\vm}_{h,\cdot}\overline{\vm}_{h,\cdot}^{\tau}$ of $var(\vm(y))$ is consistent, we must require that the average loss of variance in each slice to decrease to zero as $H$ increases, i.e., 
\begin{equation}
\begin{aligned}\label{eq:inlince:sliced}
\frac{1}{H}\sum_{h}\bar{\vm}_{h,\cdot}\bar{\vm}^{\tau}_{h,\cdot}-\frac{1}{n}\sum_{i}\vm_{i}\vm_{i}^{\tau}
=\frac{1}{H}\sum_{h}\frac{1}{c}\sum_{i}(\bar{\vm}_{h,\cdot}-\bar{\vm}_{h,i})^{2} \rightarrow 0.
\end{aligned}
\end{equation}
In Definition \ref{Definition:Sliced}, we simply choose the decreasing rate to be a power of $H$. It would be easily seen that if $\vm$ is smooth and $y$ is compactly supported then \eqref{eq:inlince:sliced} holds automatically. In this sense, for general curve $\vm$ and random variable $y$, the sliced stability is a condition on smoothness of the central curve $\vm$ and tail distribution of $\vm(y)$. This is not surprised at all, since most work on the consistency of SIR estimate requires some kind of  smoothness for the central curve and a  tail distribution  control for $\vm(y)$. 
\end{rmk}

The most popular smoothness and tail condition might be the one proposed by \cite{hsing1992asymptotic} (later used in \cite{Zhu:Miao:Peng:2006}, \cite{Zhu:Ng:1995}) in their proof of the consistency of SIR, which is explained below. For $B>0$ and $n\geq 1$, let $\Pi_{n}(B)$ be the collection of all the $n$-point partitions $-B\leq y_{(1)} \leq \cdots \leq y_{(n)} \leq B$ of $[-B,B]$.  First, they assumed that the central curve $\vm(y)$ satisfies the following smoothness condition
\begin{equation*}
\lim_{n\rightarrow \infty} \sup_{y\in \Pi_{n}(B)}n^{-1/4}\sum_{i=2}^{n}\|\vm(y_{i})-\vm(y_{i-1})\|_{2}=0, \forall B>0.
\end{equation*}
Second, they assumed that for $B_{0}>0$, there exists a non-decreasing function $\widetilde{m}(y)$ on $(B_{0},\infty)$, such that
\begin{align}
&\widetilde{m}^{4}(y)P(|Y|>y) \rightarrow 0 \mbox{ as } y\rightarrow \infty\label{tail:inline}\\
\|\vm(y)-\vm(y')\|_{2} &\leq |\widetilde{m}(y)-\widetilde{m}(y')| \mbox{ for } y, y' \in (-\infty,-B_{0})\cup(B_{0},\infty)\nonumber
\end{align}
 By changing the tail condition \eqref{tail:inline} to a slightly stronger condition $\bbE[\widetilde{m}(y)^{4}]<\infty$,  \cite{neykov2015signed} proved that the modified condition implies the sliced stability condition. Now, we are ready to state our main results.
%



\begin{thm}\label{thm:high:dim}
 Under conditions ${\bf (A1), (A2), (A3) }$ and ${\bf (A4)}$, we have
 \begin{equation}\label{eqn:sir_estimator:conditional_covariance_matrix}
\|\widehat{\bLambda}_{H}-\bLambda_{p}\|_{2}= O_{P}(\frac{1}{H^{\vartheta}}+\frac{H^{2}p}{n}+\sqrt{\frac{H^{2}p}{n}}).
 \end{equation}
 \end{thm}
The proof of the theorem is deferred to the Appendix.
As a direct consequence of Theorem \ref{thm:high:dim}, we observe that if $\rho=\lim_{n\rightarrow \infty}\frac{p}{n} =0$, we may choose $H=\log\left({n}/{p}\right)$ such that the right hand side of equation  \eqref{eqn:sir_estimator:conditional_covariance_matrix} converges to 0. Thus, Theorem~\ref{thm:high:dim} implies that $\widehat{\bLambda}_{H}$ is a consistent estimate of $\bLambda_{p}$ if $\rho=0$. 

\begin{rmk}[More on Convergence Rate]
Note that the convergence rate in \eqref{eqn:sir_estimator:conditional_covariance_matrix} depends on the choice of $H$. This may seem not very desirable at the first glance.  
Since the convergence rate of $\widehat{\bLambda}_{H}$ might be different from that of $col(\widehat{\vV}_{H})$, we may expect that the convergence rate of $col(\widehat{\vV}_{H})$ does not depend on the choice of $H$. In fact, we have
\begin{align}\label{eqn:temp:temp:rmk}
\widehat{\bLambda}_{H}-\bLambda_{p} = \left( \widehat{\bLambda}_{H}-P_{\vV}\widehat{\bLambda}_{H}P_{\vV}\right)+\left(P_{\vV}\widehat{\bLambda}_{H}P_{\vV}-\bLambda_{p}\right).
\end{align}
From the proof of Theorem 1, we can easily check that
 the first term is of convergence rate $\frac{pH^{2}}{n}+\sqrt{\frac{pH^{2}}{n}}$ and the second term is of rate $\frac{1}{H^{\vartheta}}$.  Since $P_{\vV}\widehat{\bLambda}_{H}P_{\vV}$ and $\bLambda_{p}$ share the same column space, if we are only interested in estimating $P_{\vV}$, then  the convergence rate of the second term does not matter provided that $H$ is a large enough integer, which may depend on $\vartheta$ and $\bgamma_{3}$ but does not depend on $n$ and $p$.  For such an $H$,  if $\mathcal{A}_{H}(\bgamma_{1},\bgamma_{2})$ is non-empty, Theorem \ref{thm:high:dim} and \eqref{eqn:temp:temp:rmk} hold for both categorical and continuous response variable $Y$.
\end{rmk}


\begin{thm}\label{thm:space:consistency}
Under conditions ${\bf (A1), (A2), (A3) }$, ${\bf (A4)}$ and assuming that $\rho=\lim\frac{p}{n}=0$, we have
\[
\|\widehat{\bSigma}^{-1}_{\vx}\widehat{\bLambda}_{H}-\bSigma^{-1}_{\vx}\bLambda_{p}\|_{2}\rightarrow 0 \quad \mbox{ as } n \rightarrow \infty
\]
with probability converging to one, where $\widehat{\bSigma}_{\vx}=\frac{1}{n}\sum_{i=1}^{n}\vx_{i}\vx^{\tau}_{i}$.
\end{thm}

 We define the distance $\mathcal{D}( \vV_{1}, \vV_{2})$ of two $d$-dimensional subspaces $\vV_{1}$ and $\vV_{2}$ as the operator norm (or Frobenius norm) of the difference between $P_{\vV_{1}}$ and $P_{\vV_{2}}$.
Simple linear algebra shows that if the $\widetilde{\bbeta}_{i}$'s  satisfy
$
\bSigma_{\vx}\widetilde{\bbeta}_{i}=\lambda_{i}\bLambda_{p}\widetilde{\bbeta}_{i},
$
 then 
 \[
 col(\vV)=span\Big\{\widetilde{\bbeta}_1,\cdots,\widetilde{\bbeta}_d\Big\}.
 \]
 Let $\widehat{\vV}$ be the matrix formed by the top $d$ generalized eigenvectors of $(\widehat{\bSigma}^{-1}_{\vx}, \widehat{\bLambda}_{H})$. 
Recall that the $d$-th eigenvalue of $\bLambda_{p}$ is assumed to be bounded away from $0$. Therefore Theorem \ref{thm:space:consistency} implies that $\mathcal{D}( P_{\widehat{\vV}},  P_{\vV})\to 0$ when $\rho=0$.


We have already shown that the SIR procedure provides us with a consistent estimate of the sufficient dimension reduction space when $\rho = 0$ under mild conditions. 
It is then natural to ask: is this condition necessary? Our next theorem gives the answer.

\begin{thm} \label{thm:main1.1} Under conditions ${\bf (A1), \ (A2), \ (A4)}$ and assuming that $\vx \sim N(0,\vI_{p})$ for the single index model
\[
y=f(\bbeta^{\tau}\vx,\epsilon),
\] 
we have:
\begin{itemize}
\item[(i)] When $\rho=\lim\frac{p}{n}\in(0,\infty)$, $\|\widehat{\bLambda}_{H}-\bLambda_{p}\|_{2}$, as a function of $\rho$, is dominated by $\sqrt{\rho}\vee\rho$ when $H, n \rightarrow \infty$; 
\item[(ii)] Let $\widehat{\bbeta}$ be the principal eigenvector of the SIR estimator $\widehat{\bLambda}_{H}$. If $\rho=\lim\frac{p}{n}>0$, then there exists a positive constant $c(\rho)>0$ such that 
\[
\liminf_{n\rightarrow \infty} \bbE\angle(\bbeta,\widehat{\bbeta})>c(\rho)
\]
with probability converges to one.
\end{itemize}
\end{thm}

\begin{figure}[!h]
\centering
\includegraphics[height=30mm,width=60mm]{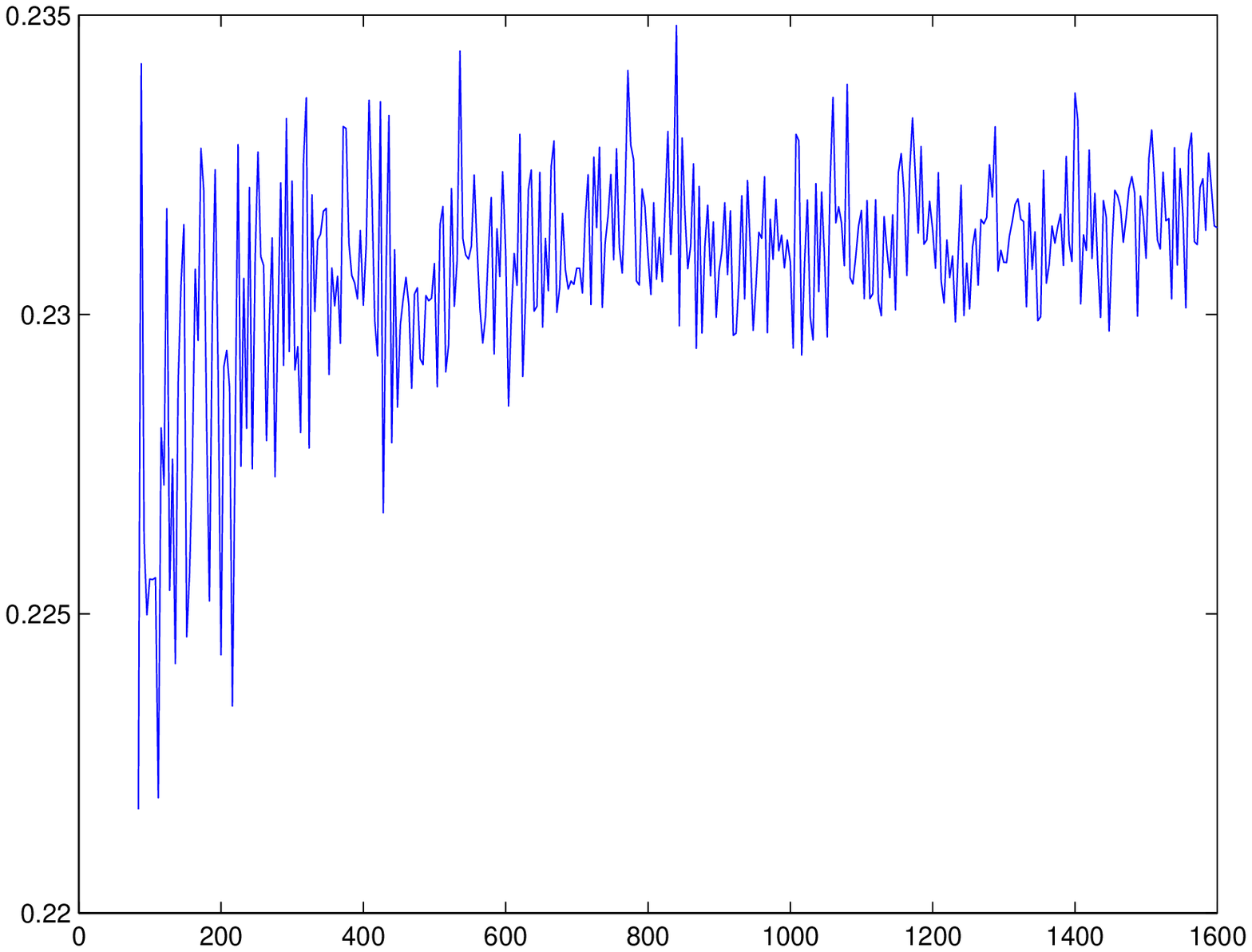}
\includegraphics[height=30mm,width=60mm]{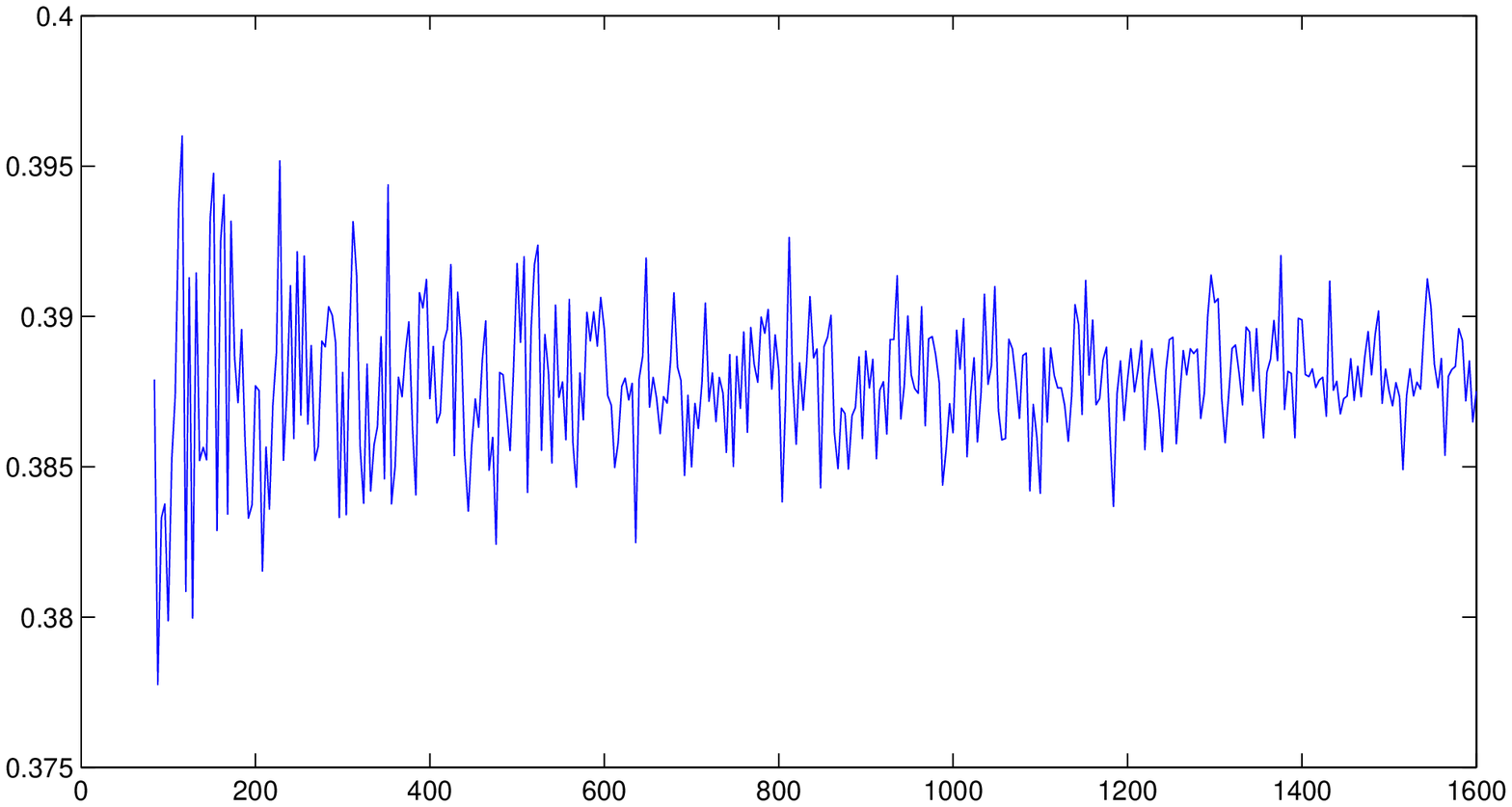}
\includegraphics[height=30mm,width=60mm]{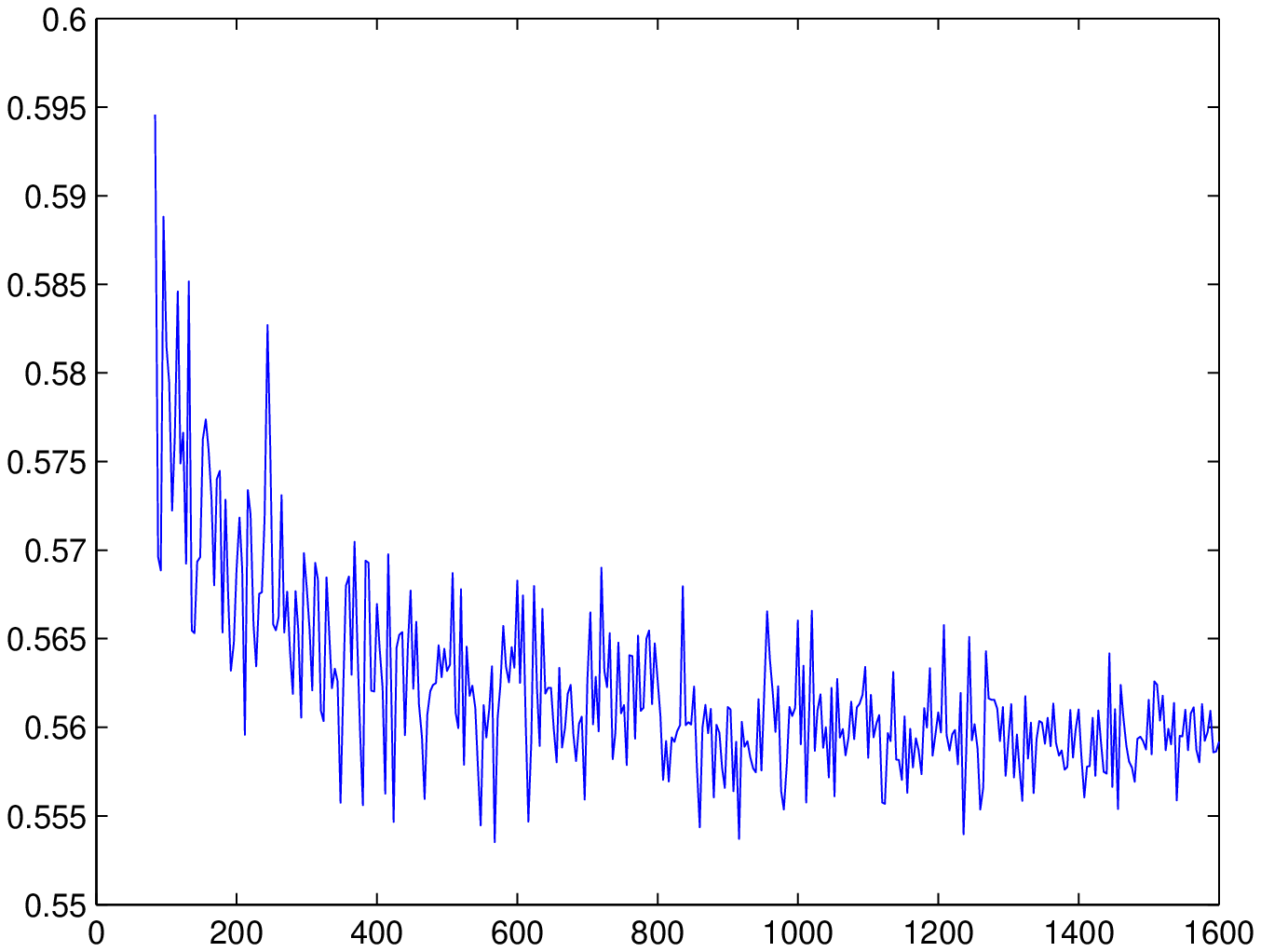}
\includegraphics[height=30mm,width=60mm]{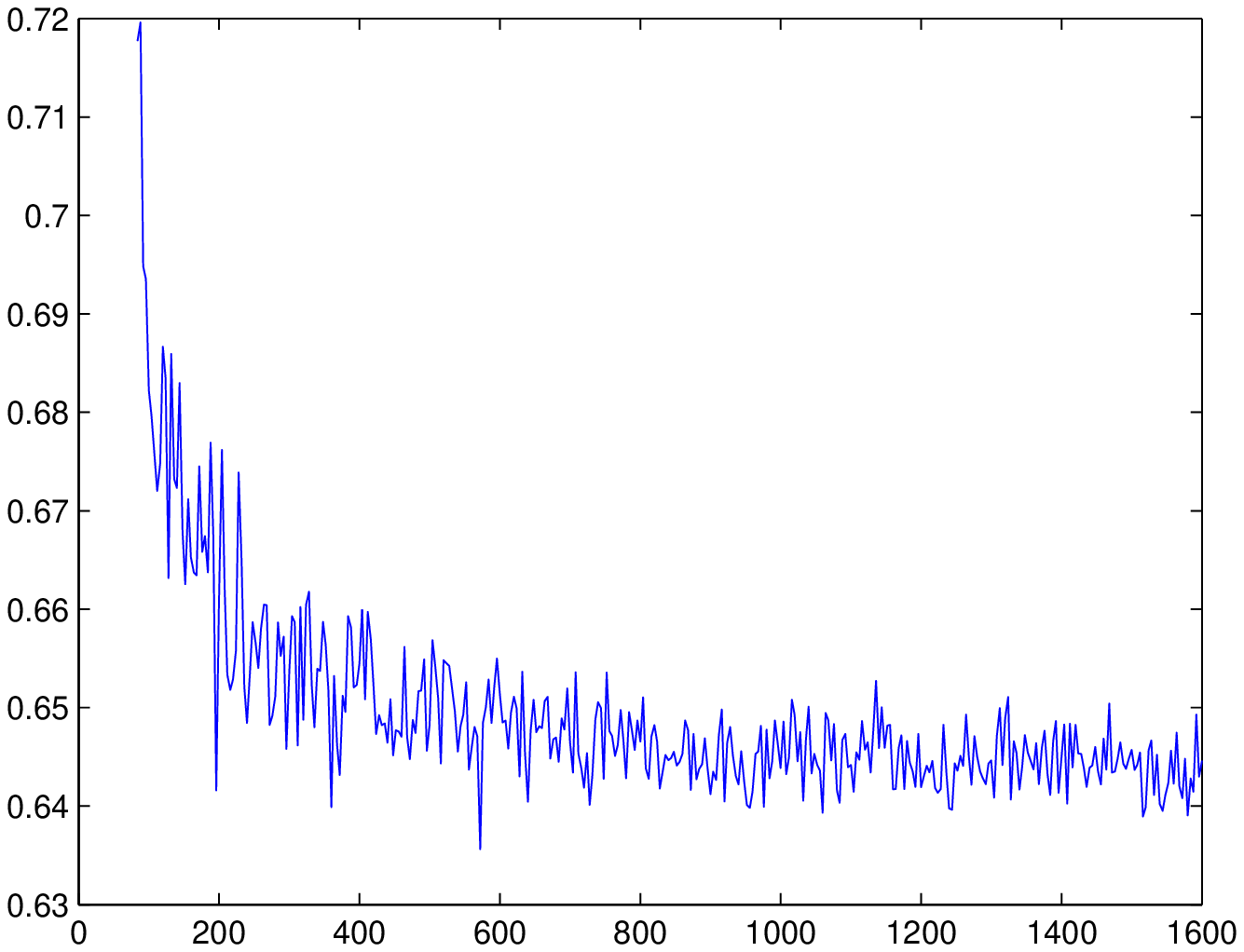}
\includegraphics[height=30mm,width=60mm]{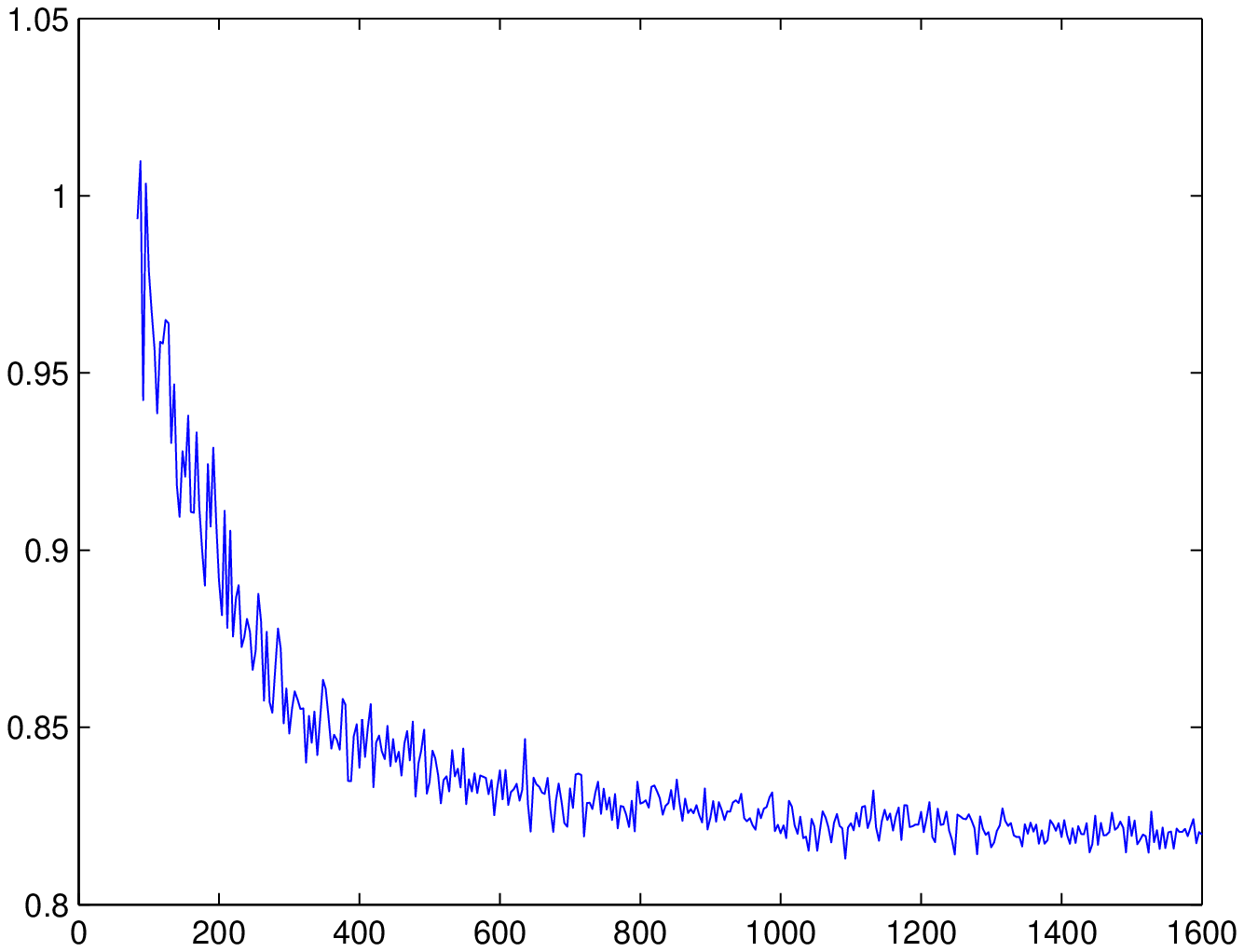}
\includegraphics[height=30mm,width=60mm]{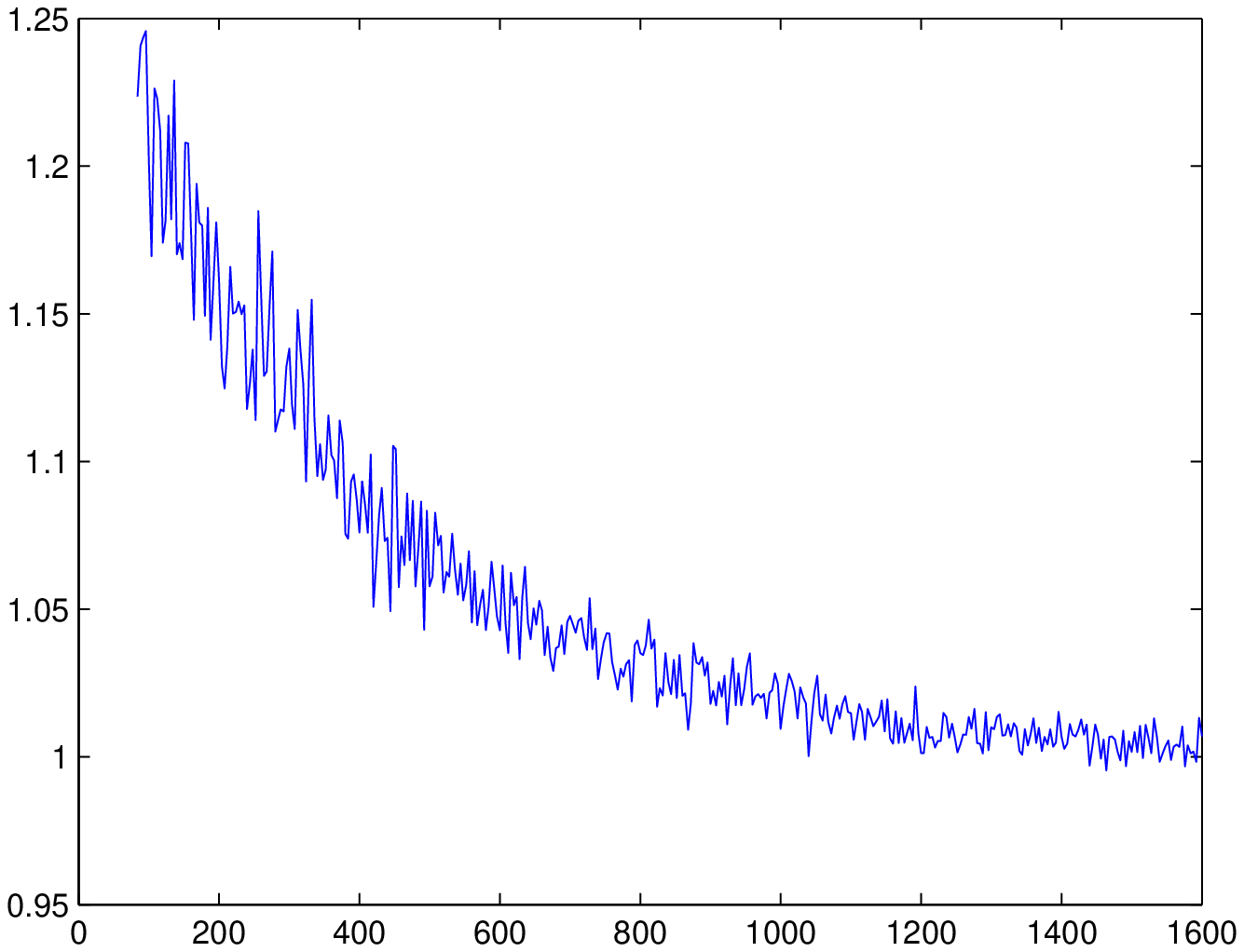}
\caption{Numerical approximations of $\bbE\angle(\hat{\bbeta}, \bbeta)$ for model \eqref{linear-model}  as a function of  dimension $p$ for 
$\rho=.1$, .3, .7, 1, 2, and 4, respectively (up left, up right, middle left, middle right, lower left, lower right), where $\widehat{\bbeta}$ is estimated by SIR.}\label{fig:9}
\end{figure}

\begin{figure}[!h]
\centering
\includegraphics[height=50mm, width=100mm]{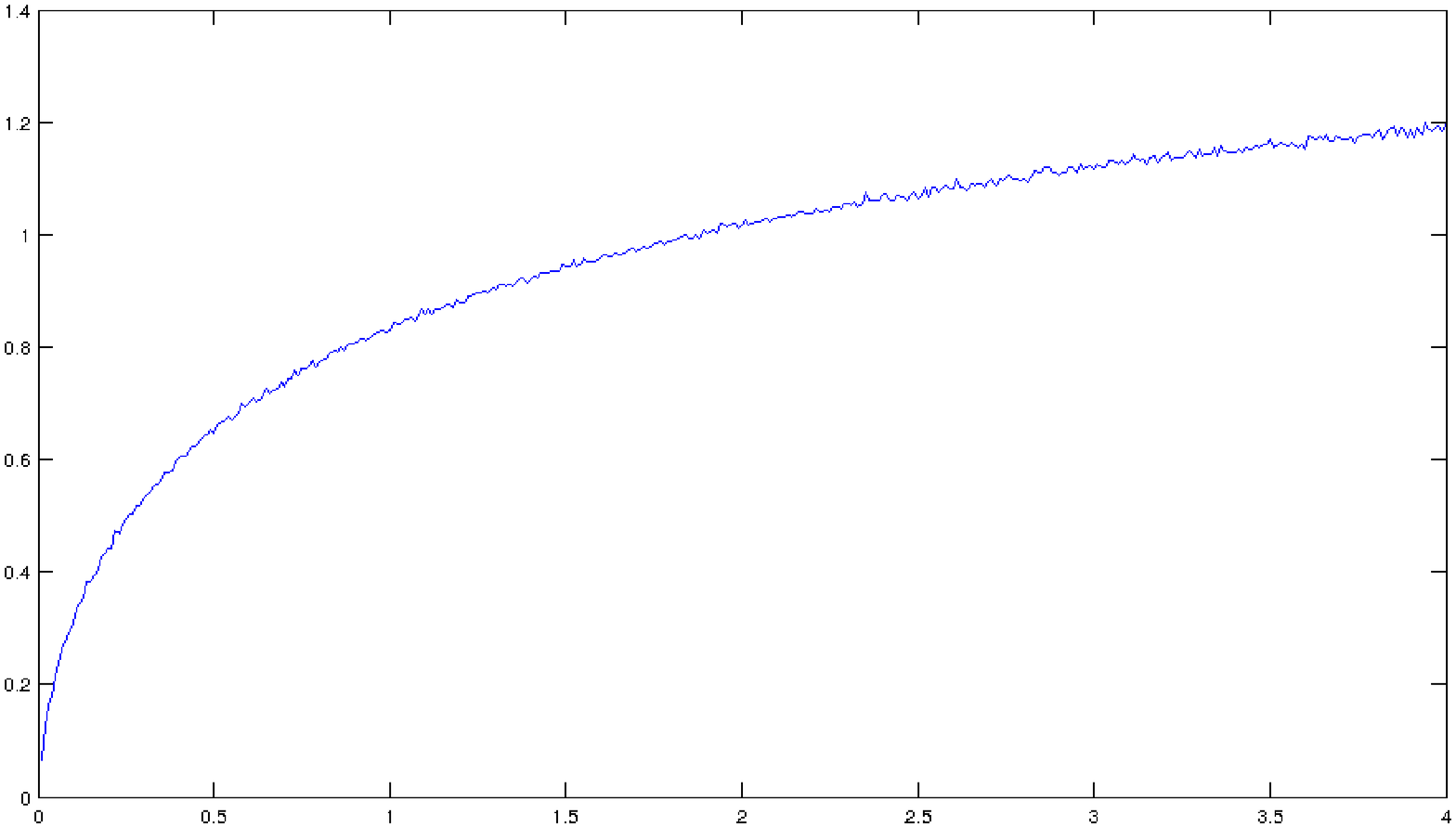}
\caption{The relationship of $\bbE\angle(\bbeta, \widehat{\bbeta})$ and the ratio $p/n$ where $\widehat{\bbeta}$ is estimated by SIR.}\label{fig:1}
\end{figure}

We illustrate this result via a numerical study of the  linear model
\begin{equation}\label{linear-model}
y=\vx^{\tau}\bbeta+\epsilon \mbox{ where } \bbeta^{\tau}=(1,0,\cdots,0), \vx \sim N(0,\vI_{p}), \epsilon \sim N(0,1).
\end{equation}
Figure \ref{fig:9} shows how  $\bbE\angle(\bbeta, \widehat{\bbeta})$ is related to the dimension $p$ for fixed ratio $\rho=\frac{p}{n}$ (taking values in  $\{.1,.3,.7,1,2,4\}$),  where  $\bbeta$ is estimated by the SIR with the slice number $H=10$. For each $p$,  $\bbE\angle(\bbeta, \widehat{\bbeta})$ is calculated based on 100 iterations. It is seen that this expected angle converges to a positive number when the ratio $\rho$ is non-zero.
In Figure \ref{fig:1}, we have plotted the $\bbE\angle(\bbeta, \widehat{\bbeta})$ against the ratio $\rho=\frac{p}{n}$, varying between 0.01 and 4 with an increment of 0.01. 
The sample size $n$ is 200 and the slice number $H$ is 10. 
It is seen that the expected angle decreases to zero as $\rho$ approaches zero, and  increases monotonically when $\rho$ increases.

Results in this section have shown that there is a phase transition phenomenon  of the SIR procedure. 
That is, the estimate of the dimension reduction space is consistent if and only if the ratio $\rho=\lim\frac{p}{n}=0$. 
This provides a theoretical justification of  imposing additional structure assumption such as sparsity  in high dimension.

\section{SIR in ultra-high dimension.}\label{sec:high_dimensional}
As we have shown in Section \ref{sec:sir}, the SIR estimator fails to be consistent if $\rho=\lim\frac{p}{n}\neq 0$. 
Hence, when $p\gg n$,  some structural assumptions are necessary for getting a consistent estimate of the central space.  
In this paper, we assume that  both  the loadings of all the directions $\bbeta_j$'s  and  the covariance matrix $\bSigma_{\vx}$ are sparse. 
Other structural assumptions will be studied in our future work. For $\bbeta_{i}$'s, we impose the following prevalent sparsity condition.
\begin{itemize}
\item {\bf (A5)} $s=|\mathcal{S}|\ll p$ where $\mathcal{S}=\Big\{~i~| \textrm{ $\bbeta_{j}(i)\neq 0$ for some $j$, $1\leq j\leq d$ }\Big\}$ and $|\mathcal{S}|$ is the number of elements in the set $\mathcal{S}$.
\end{itemize} 
For  $\bSigma_{\vx}$, the following class of covariance matrices has been introduced in \cite{bickel2008covariance} (see also \cite{Cai:Zhang:Zhou:2010}).
\begin{align*}
\mathcal{U}(\epsilon_{0},\alpha,C)
=\Big\{~\bSigma_{\vx}: 
\max_{j}\sum_{i}\{|\sigma_{i,j}|: |i-j|>l\} \leq Cl^{-\alpha} \mbox{ for all } l>0, \\
 \mbox{ and } 0< \epsilon_{0}\leq \lambda_{min}(\bSigma_{\vx})\leq \lambda_{max}(\bSigma_{\vx})\leq \frac{1}{\epsilon_{0}}
~\Big\}.
\end{align*}
In this paper, to simplify the notations and arguments, we choose a slightly stronger condition.
\begin{itemize}
\item {\bf (A6)} $\bSigma_{\vx} \in \mathcal{U}(\epsilon_{0},\alpha,C)$ and $\max_{1\le i\le p}r_i$ is bounded where $r_{i}$ is the number of non-zero elements in the i-th row of $\bSigma_{\vx}$.
\end{itemize}

Let $\mathcal{T}=\big\{\Scale[.9]{~k~|~ var\left(\bbE[\vx(k)|y]\right)\neq 0~}\big\}$. 
If $k \in \mathcal{T}$, there exists $\boldeta \in col(\bLambda)$ such that $\boldeta(k)\neq 0$. 
Since we have \eqref{eta:beta}:
\begin{align*}
\bSigma_{\vx}col(\vV)=col(\bLambda),
\end{align*}
 there exits a $\bbeta \in col(\vV)$ such that 
$\boldeta=\Sigma_{\vx}\bbeta$. 
Thus if $k \in \mathcal{T}$,  then $k \in supp(\Sigma_{\vx}\bbeta)$ for some $\bbeta \in col(\vV)$. In particular, with the above sparsity assumptions {\bf (A5)}  and {\bf (A6)} ,  we have $|\mathcal{T}|\leq s\max_{1\leq i \leq p}r_{i}=O(s)$.\footnote{We could introduce $\xi=\max_{1\leq i \leq p}r_{i}$, then $|\mathcal{T}|\leq s\xi$. The arguments below still work, except we might need $s\xi =o(p)$.} Note that our goal here is to recover the column space $col(V)$ rather than $\mathcal{S}$. Indeed, we are not able to consistently recover $\mathcal{S}$ unless for the trivial case. The key for recovering  $cov(V)$ is to consistently recovering the set $\mathcal{T}$.

At the population level,  $var(\bbE(\vx(k))|y)$ can separate $\mathcal{T}$ from $\mathcal{T}^c$.  
When there are only finite samples, we use 
\begin{align}\label{statistics:screening}
var_{H}(\vx(k)) = \frac{1}{H} \sum_{h=1}^H  \bar{\vx}_{h,\cdot}(k)^2
\end{align}
as an estimate of $var(\bbE(\vx(k))|y)$.
These are the diagonal elements of the matrix $\widehat{\bLambda}_{H}$. 
Note that these  quantities depend on  the sliced sample means, which are neither independent nor identically distributed. Thus, the usual concentration inequalities for $\chi^{2}$ are no longer applicable.  We need extra efforts to get the concentration inequalities; this concentration result is one of the main technical contributions of this article, and can be further generalized. 

\begin{rmk}
  The link function $f(\ )$ is not involved explicitly in the definition of $var_{H}(\vx(k))$, and only the order statistics of the response is required. This nonparametric characteristic of the method is of particular interest to us and will be further investigated in a future research. Screening statistics inspired by the sliced inverse regression idea have been proposed in various formats, such as those in \cite{jiang2013sliced}, \cite{zhu2011model}  and \cite{cui2014model}. 
\end{rmk}

With the quantities  $var_{H}(\bbE[\vx(k)|y])$, we define the inclusion set $\mathcal{I}_{p}(t)$ and the exclusion set $\mathcal{E}_{p}(t)$ below, which depend on a thresholding value t:
\begin{align*} 
\mathcal{I}_{p}(t)=\Big\{~k~| var_{H}(\vx(k)) > t ~\Big\} \mbox{ and } \mathcal{E}_{p}(t)=\Big\{~k~| var_{H}(\vx(k)) \leq t ~\Big\}.
\end{align*}
Note that $\mathcal{I}_{p}(t)$ can be viewed as an estimate of $\mathcal{T}$ and is thus also denoted by $\widehat{\mathcal{T}}$. 
After reducing the dimension to a level such as $p/n$ is sufficiently small, 
 the SIR estimator $\Scale[0.8]{\widehat{\bLambda}^{\widehat{\mathcal{T}}, \widehat{\mathcal{T}}}}$ is a consistent estimate of $\Scale[0.8]{\bLambda^{\mathcal{T},\mathcal{T}}}$. 
Let $\Scale[0.8]{\widehat{\vV}^{\widehat{\mathcal{T}}}}$ be the matrix formed by the  top $d$ eigenvectors of $\Scale[0.8]{\widehat{\bLambda}^{\widehat{\mathcal{T}},\widehat{\mathcal{T}}}}$.  
We then use $\Scale[0.8]{\widehat{\bSigma}^{-1}_{\vx}}col(\Scale[0.8]{e(\widehat{\vV}^{\widehat{\mathcal{T}}})})$ to estimate the central space $col(\vV)$, where $\Scale[0.8]{\widehat{\bSigma}^{-1}_{\vx}}$ is a consistent estimate of $\bSigma_{\vx}$.  Estimating the covariance matrix and precision matrix in high dimension is a challenging problem by itself and is not a main focus of this article. We employ the methods of \cite{bickel2008covariance} to solve it. In summary, we propose the following {\bf D}iagonal {\bf T}hresholding screening {\bf SIR} (DT-SIR)  algorithm:

\begin{algorithm}[H]
\label{CHalgorithm}
\begin{algorithmic}
\caption{DT-SIR}

\vspace{3mm}
\item[1.] Calculate $var_{H}(\vx(k))$ according \eqref{statistics:screening} for $k=1,2,\cdots, p$;

\vspace{3mm}
\item[2.] Let $\widehat{\mathcal{T}}=\Big\{~k~| ~var_{H}(\vx(k)) > t ~\Big\}$ for an appropriate t ;

\vspace{3mm}
\item[3.] Let  $\Scale[0.8]{\widehat{\bLambda}_{H}^{\widehat{\mathcal{T}},\widehat{\mathcal{T}}}}$ be the SIR estimator of the conditional covariance matrix for the data  $(y,\Scale[0.8]{\vx^{-,\widehat{\mathcal{T}}}})$ according to equation \eqref{eqn:lambda}; 

\vspace{3mm}
\item[4.]  Let $\Scale[0.8]{\widehat{\vV}^{\widehat{\mathcal{T}}}}$ be the matrix formed by the   top $d$ eigenvectors of $\Scale[0.8]{\widehat{\bLambda}^{\widehat{\mathcal{T}},\widehat{\mathcal{T}}}}$;

\vspace{3mm}
\item[5.] $\Scale[0.8]{\widehat{\bSigma}^{-1}_{\vx}}col\left(e\left(\Scale[0.6]{\widehat{\vV}^{\widehat{\mathcal{T}}}}\right)\right)$ is the estimate of $col(\vV)$ 

\end{algorithmic}
\end{algorithm}

A practical way to choose an  appropriate $t$ in step 2 will be presented in Section \ref{sec:simulation} .
To ensure theoretical properties, we need an  assumption on the signal strength:
\begin{itemize}
\item ({\bf S1}) $\exists$ $C>0$ and $\omega>0$ such that $var(\bbE[\vx(k)|y])>{C}{s^{-\omega}}$
when $\bbE[\vx(k)|y]$ is not a constant.
\end{itemize}
\begin{thm}\label{thm:screen:consistent}
Under conditions  {\bf (A1) -- (A6)} and  {\bf (S1)}, and let  $t={a}{s^{-\omega}}$ for some  constant $a>0$ such that $t< \frac{1}{2}var(m(y,k), \ \forall k \in \mathcal{T}$, we have 
\begin{itemize}

\item[i)] $\mathcal{T}^{c} \subset \mathcal{E}_{p}$ holds with probability at least
\end{itemize}
\begin{equation}\label{ex_prob}
\begin{aligned}
 1-C_{1}\exp\left(-C_{2}\frac{n}{H^{2}s^{\omega}}+C_{3}\log(H)+\log(p-s)\right);
\end{aligned}
\end{equation}
\begin{itemize}
\item[ii)] $\mathcal{T} \subset \mathcal{I}_{p}$ holds with probability at least
\end{itemize}
\begin{equation}\label{in_prob}
\begin{aligned}
1-C_{4}\exp\left(-C_{5}\frac{n}{H^{2}s^{\omega}}+C_{6}\log(H)+\log(s)\right),
\end{aligned}
\end{equation}
for some positive constants $C_{1}, \cdots, C_{6}$.
\end{thm}
This theorem has a simple implication. If 
$
\frac{n}{s^{\omega}} \succ \log(p)+\log(s)
$,
we may choose $H=\log(\frac{n}{s^{\omega}\log(p)})$, so that 
\[
\frac{n}{H^{2}s^{\omega}} \succ \log(p)+\log(H)+\log(s).
\]
Thus , we know $\mathcal{T}=\mathcal{I}_{p}$ with probability converging to one. 
Next, we have results for the consistency of   DT-SIR. 

\begin{thm} \label{thm:conditional_variance:high_dimensional}
Under the same assumptions and choosing the same $t$ as Theorem \ref{thm:screen:consistent}, if $
\frac{n}{s^{\omega}} \succ \log(p)+\log(s)
$,  we have   
\[
\|e(\Scale[0.8]{\widehat{\bLambda}_{H}^{\widehat{\mathcal{T}},\widehat{\mathcal{T}}}})-\bLambda_p\|_{2} \rightarrow 0  \mbox{ as } n \rightarrow \infty
\]
with probability converging to 
one, where  $\widehat{\mathcal{T}}=\mathcal{I}(t)$ and $H=\log(\frac{n}{s^{\omega}\log(p)})$.
\end{thm}

\begin{thm}\label{thm:operator_norem:high_dimensional} Let $\widehat{\bSigma}_{\vx}$ be the estimator of co-variance matrix from \cite{bickel2008covariance}. Under the same assumptions of Theorem \ref{thm:conditional_variance:high_dimensional}, we have 
\[
\|\widehat{\bSigma}^{-1}_{\vx}e(\Scale[0.8]{\widehat{\bLambda}_{H}^{\widehat{\mathcal{T}},\widehat{\mathcal{T}}}})-\bSigma^{-1}_{\vx}\bLambda_p\|_{2} \rightarrow 0 \mbox{ as } n \rightarrow \infty
\]
with probability converging to one.
\end{thm}

\section{Simulation Studies.}\label{sec:simulation}

We consider the following settings in generating the design matrix $\vx$ and the response $y$.
In Settings I-III, each row of $\vx$ is independently sampled from $N(\vzero,\vI)$.
\begin{itemize}
\item {\bf Setting I.}  $y_i=\sin(x_{i1} + x_{i2})+ \exp(x_{i3}+x_{i4}) + 0.5* \epsilon_i$, where $\epsilon_i\iid N(0,1)$;
\item {\bf Setting II.} $y_i = \sum_{j=1}^7x_{ij}  * \exp( x_{i8}+x_{i9}) + 0.5 * \epsilon_i$ where $\epsilon_i\iid N(0,1)$;
\item {\bf Setting III.} $y_i=\sum_{j=1}^{10}x_{ij} *\exp (\sum_{i=11}^{20}x_{ij}) + \epsilon_i$ where $\epsilon_i\iid N(0,1)$;
\end{itemize}
In Settings IV to VI, each row of $\vx$ is independently sampled from $N(\vzero,\bSigma)$.
\begin{itemize}
\item {\bf Setting IV.} $y_i=(x_{i1}+x_{i2}+x_{i3})^3/2+ 0.5* \epsilon_i$, where $\epsilon_i\iid N(0,1)$ and $\bSigma=(\sigma_{ij})$ is tri-diagonal with $\sigma_{ii}=1$, $\sigma_{i,i+1}=\sigma_{i+1,i}=\rho$ and $\sigma_{i,i+2}=\sigma_{i+2,i}=\rho^2$;
\item {\bf Setting V.} $y_i = \sum_{j=1}^7x_{ij}  * \exp( x_{i8}+x_{i9}) + \epsilon_i$, where $\epsilon_i\iid N(0,1)$, and $\bSigma= \vB \otimes \vI_{p/10}$ with $\vB=(b_{ij})_{1\le i\le 10, 1\le j\le 10}$ given as $b_{ij}=\rho^{|i-j|}$;
\item {\bf Setting VI.} Assume the same setting as in Setting V except that $\bSigma=(\sigma_{ij})$ is tri-diagonal with $\sigma_{ii}=1, \sigma_{i,i+1}=\sigma_{i+1,i}=\rho$ and $\sigma_{i,i+2}=\sigma_{i+2,i}=\rho^2$.
\item {\bf Setting VII.} Assume the same setting as in Setting V except that $\bSigma=(\sigma_{ij})$ is given as $\sigma_{ij}=\rho^{|i-j|}$.

\end{itemize}

  DT-SIR first screens all the predictors according to the statistic $var_{H,c}(\vx(k))$, which requires a tuning parameter $t$. We chose $t$ by using an auxiliary variable method based on an idea first proposed by \cite{Luo:Stefanski:Boos:2006} and extended by \cite{Wu:Boos:Stefanski:2007} and \cite{Zhu:Li:Li:Zhu:2011}. In our setting, for a given sample $(y_i, \vx_i)$, we generate $\vz_i\sim N(\vzero, \vI_{p'})$ where $p'$ is sufficiently large and chosen as $p$ in our simulation studies. It is known that $\vy$ and $\vz$ are independent. The threshold $t$ can be chosen as 
\[
\hat{t}=\max_{1\le k\le p'}\{ var_{H,c}(\vz(k)). \}
\]
In DT-SIR, when $n>1000$, $H$ is chosen as 20; when $n\le 1000$, $H$ is chosen as 10 in the  screening step and 20 in the SIR step.




We also consider the following alternative methods in the screening step: Sure Independent Ranking and Screening (SIRS) in \cite{Zhu:Li:Li:Zhu:2011}, SIR for variable selection via Inverse modeling (SIRI) in \cite{jiang2013sliced}, and trace pursuit in \cite{yu2015trace}. As a comparison, we also considered two screening methods that are not based on the sliced regression: Distance correlation in \cite{szekely2007measuring} and SURE independence \cite{fan2008sure}. For SIRS, the threshold is chosen according to the auxiliary statistic (2.9) of \cite{Zhu:Li:Li:Zhu:2011}. For SIRI, the predictors are chosen according to 10-fold cross validation. The threshold values $\bar{c}^{SIR}$ and $\underline{c}^{SIR}$ are chosen as the 10-th and 5-th quantile of a weighted $\chi^2$ distribution given in Theorem 3.1 of \cite{yu2015trace}. In both SURE and DC screening, the top $\lfloor \gamma n\rfloor$ where $\gamma=0.01$ are kept for subsequent analyses.

After the screening step, similar to DT-SIR, we then applied the SIR  algorithm (steps 3-5 of DT-SIR) to estimate $col(\vV)$. These alternative  methods are denoted as SIRS-SIR, SIRI-SIR, SURE-SIR, DC-SIR, and TP-SIR, respectively,  in the following discussions. 
Another method that we compared with is the sparse SIR, abbreviated as SpSIR, proposed in \cite{li2007sparse}.
After obtaining an estimator $col(\widehat{\vV})$, we calculate $\mathcal{D}(P_{col(\Scale[0.6]{\widehat{\vV}})}, P_{col(\vV)})$ as a measure of the estimation error. 
We replicate this step 100 times, and calculate the average distance for the estimation result from each method and report these numbers in Table \ref{table:sim:1}-\ref{table:sim:3}. 
For each setting, the average distance of the optimal method is highlighted using bold fonts.
We further run a two-sample T-test to test if the actual estimation error  of each method is significantly different from that of the best method for that example at 1\% level of significance.



Under all settings, the average distance obtained by  DT-SIR was much smaller than that obtained by SpSIR and SURE-SIR. The p-values for comparing DT-SIR and SpSIR/SURE-SIR are all significant at the 0.01 level. 
When $p\ge n$, the sparse SIR completely failed because the average distance of the estimated space to the true space is $\sqrt{2d}$, indicating that the space estimated by sparse SIR is orthogonal to the true space spanned by $\beta$. 

Under settings II-IV,  DT-SIR  performed either the best or not significantly worse than the best method. For all  other cases,  DT-SIR performed the best except for a few cases: Setting I when $n=500, p=1000$, setting V when $n=500, p=6000$, setting VI when $n=500, p=6000$, and setting VII when $n=1000, p=1000$. 

When $p=6000, n=500$, both DT-SIR and SIRI-SIR were the winners. Under Setting III, DT-SIR performed better than SIRI-SIR; under settings V and VI, SIRI-SIR performed better than DT-SIR; under other settings, these two methods were comparable.

\begin{table}[!h]
\caption{The average distance of the space estimated by each of the 7 methods tested to the true space $col(\vV)$ under various settings with $p=1000$. The boldfaced number in each row represents the best result for that simulation scenario, and the ``*'' in cells represents that the p-value of the two-sample T-test comparing the estimation error of the corresponding method with that of the best method  is less than 0.01.
}\label{table:sim:1}
\begin{tabular}{|c|c|c|c|c|c|c|c|c|}
\hline
 & n & DT-SIR & SIRI-SIR & SIRS-SIR & SpSIR & SURE-SIR & DC-SIR & TP-SIR\\
\hline
\multirow{4}{*}{I}  & 500    & 0.655(*) &0.751(*) &{\bf 0.492}&2(*) &1.39(*) &0.731(*) &1.18(*) \\
 & 1000 &   {\bf 0.3}  &  0.431(*) &0.309&2(*)  &1.29(*) &0.632(*) &0.94(*) \\
 & 2000 & {\bf 0.221}  &  0.341(*) &0.226&1.58(*)   &1.04(*) &0.655(*) &0.784(*) \\
 & 3000 & 0.167  &  0.245(*) &{\bf 0.149} &1.48(*)   &0.816(*) &0.641(*) &0.713(*)  \\
\hline
\hline
\multirow{4}{*}{II} & 500  &  0.383  &  0.396& {\bf 0.371} &2(*)  &1.64(*) &1.08(*) &0.389\\
 & 1000 & 0.235  & {\bf  0.227} &0.256 &2(*)   &1.36(*) &0.266(*) &0.318(*)   \\
 & 2000 & 0.161  &  {\bf 0.157} &0.189(*) &1.25(*)  &1.25(*) &0.387(*) &0.264(*) \\
 & 3000 & 0.134  &  {\bf 0.129} &0.153(*) &0.975(*) &1.12(*) &0.404(*) &0.23(*)  \\

\hline
\hline
\multirow{4}{*}{III} & 500  & 1.15  &  1.48(*) &1.38(*) &2(*)&1.97(*) &1.85(*) &{\bf 1.13} \\
 & 1000 &  {\bf 0.426}  &  0.974(*) &0.596(*) &2(*)  &1.94(*) &1.57(*) &0.429 \\
 & 2000 & {\bf 0.263}  &  0.403(*) &0.29(*) &1.33(*) &1.89(*) &0.996(*) &0.338(*) \\
 & 3000 & {\bf 0.214 }  &  0.297 &0.238(*) &1.06(*)  &1.82(*) &0.475(*) &0.299(*)  \\
\hline
\hline
\multirow{4}{*}{IV} &  500  &   0.263  &  {\bf 0.257} &0.333 &1.41(*)&0.335(*) &0.334(*) &0.332(*) \\
 & 1000 & {\bf 0.219}  &  0.447(*) &0.25 &1.41(*) &0.436(*) &0.459(*) &0.469(*)  \\
 & 2000 & {\bf 0.161}  &  0.4(*) &0.196(*) &0.42(*)   &0.442(*) &0.469(*) &0.452(*) \\
 & 3000 & {\bf 0.134}  &  0.377(*) &0.177(*) &0.297(*) &0.43(*) &0.458(*) &0.438(*)  \\
\hline
\hline
\multirow{4}{*}{V} & 500  &   0.546  &  {\bf 0.529} &0.562&2(*) &1.62(*) &1.24(*) &1.09(*)  \\
 & 1000 & 0.401  &  0.463(*) &0.514(*) &2(*)  &1.15(*) &{\bf 0.367}&0.615(*)  \\
 & 2000 & {\bf 0.288}  &  0.418(*) &0.341(*) &1.51(*) &0.926(*) &0.569(*) &0.54(*) \\
 & 3000 & {\bf 0.249}  &  0.399(*) &0.284(*) &1.24(*)  &0.691(*) &0.597(*) &0.511(*)  \\
\hline
\hline
\multirow{4}{*}{VI} & 500  &   0.568  & {\bf 0.535 } &0.566&2(*)  &1.64(*) &1.24(*) &1.08(*)\\
 & 1000 & 0.427  &  0.524(*) &0.548(*) &2(*)  &1.22(*) &{\bf 0.39}&0.641(*)   \\
 & 2000 & {\bf 0.311}  &  0.469(*) &0.351(*) &1.51(*)   &0.927(*) &0.598(*) &0.583(*) \\
 & 3000 &  {\bf 0.265}  &  0.456(*) &0.307(*) &1.25(*)  &0.807(*) &0.622(*) &0.56(*) \\
\hline
\hline
\multirow{4}{*}{VII} & 500    & 0.556& {\bf 0.534} &0.585(*) &2(*) &1.66(*) &1.26(*) &1.11(*)\\
 & 1000  & 0.436(*) &0.528(*) &0.545(*) &2(*) &1.22(*) & {\bf 0.39} &0.643(*) \\
 & 2000   & {\bf 0.303} &0.465(*) &0.358(*) &1.51(*) &0.747(*) &0.589(*) &0.579(*) \\
 & 3000   & {\bf 0.258} &0.468(*) &0.319(*) &1.25(*) &0.698(*) &0.63(*) &0.558(*)\\
\hline

\end{tabular}
\end{table}

\begin{table}[!h]
\caption{The average distance of the space estimated by each of the 7 methods we tested to the true space $col(\vV)$ under various settings with $n=2000$.
}\label{table:sim:2}
\begin{tabular}{|c|c|c|c|c|c|c|c|c|}
\hline
 & p & DT-SIR & SIRI-SIR & SIRS-SIR & SpSIR & SURE-SIR & DC-SIR & TP-SIR\\
\hline
\hline
\multirow{4}{*}{I}  & 500   & 0.213  &  0.312(*) &{\bf 0.206}&1.44(*) &0.903(*) &0.629(*) &0.772(*) \\
& 1000 & {\bf 0.221}  &  0.341(*) &0.226&1.58(*) &1.04(*) &0.655(*) &0.784(*)  \\
& 2000 & 0.241  &  0.29 &{\bf 0.214}&2(*) &1.07(*) &0.677(*) &0.793(*) \\
& 3000 & 0.23  &  0.278 & {\bf 0.218} &2(*)  &1.17(*) &0.683(*) &0.797(*)  \\
\hline
\hline
\multirow{4}{*}{II}  & 500   &  0.163  &  {\bf 0.16} &0.19(*) &0.83(*) &1.22(*) &0.369(*) &0.26(*)  \\
    & 1000 & 0.161  &  {\bf 0.157} &0.189(*) &1.25(*)  &1.25(*) &0.387(*) &0.264(*) \\
 & 2000 &  0.172  &  {\bf 0.159} &0.196(*) &2(*) &1.23(*) &0.404(*) &0.259(*)  \\
 & 3000 &  0.164  &  {\bf 0.158} &0.199(*) &2(*) &1.3(*) &0.414(*) &0.261(*)   \\

\hline
\hline
\multirow{4}{*}{III}  & 500  & {\bf 0.272}  &  0.353 &0.29(*) &0.916(*) &1.84(*) &0.846(*) &0.341(*)  \\
   &  1000 & {\bf 0.263}  &  0.403(*) &0.29(*) &1.33(*) &1.89(*) &0.996(*) &0.338(*) \\
 & 2000 &  {\bf 0.262}  &  0.368 &0.285(*) &2(*)  &1.92(*) &0.98(*) &0.339(*)  \\
   & 3000 & {\bf 0.269}  &  0.344 &0.291(*) &2(*)  &1.93(*) &1.09(*) &0.339(*)\\
\hline
\hline
\multirow{4}{*}{IV} & 500  &  {\bf 0.145}  &  0.409(*) &0.182(*) &0.248(*) &0.406(*) &0.433(*) &0.438(*)  \\
 &  1000 &  {\bf 0.161}  &  0.4(*) &0.196(*) &0.42(*)  &0.442(*) &0.469(*) &0.452(*)  \\
 &  2000 & {\bf 0.16}  &  0.395(*) &0.198(*) &1.41(*) &0.472(*) &0.506(*) &0.447(*)  \\
 &  3000 &  {\bf 0.15}  &  0.395(*) &0.216(*) &1.41(*) &0.49(*) &0.527(*) &0.447(*) \\
\hline
\hline
\multirow{4}{*}{V} &500   &   {\bf 0.272}  &  0.434(*) &0.353(*) &1.09(*)&0.876(*) &0.547(*) &0.539(*)  \\
 & 1000 & {\bf 0.288}  &  0.418(*) &0.341(*) &1.51(*)  &0.926(*) &0.569(*) &0.54(*)   \\
 & 2000 & {\bf  0.289}  &  0.418(*) &0.351(*) &2(*)  &0.868(*) &0.596(*) &0.537(*) \\
 & 3000 & {\bf  0.3}  &  0.417(*) &0.372(*) &2(*)  &0.968(*) &0.605(*) &0.544(*) \\
\hline
\hline
\multirow{4}{*}{VI} & 500   &  {\bf 0.307}  &  0.479(*) &0.368(*) &1.1(*) &0.858(*) &0.566(*) &0.583(*) \\
 &  1000 &{\bf 0.311}  &  0.469(*) &0.351(*) &1.51(*) &0.927(*) &0.598(*) &0.583(*)\\
 &2000 & {\bf 0.309}  &  0.461(*) &0.399(*) &2(*)&1.08(*) &0.617(*) &0.585(*) \\
 &3000 &  {\bf 0.31}  &  0.46(*) &0.408(*) &2(*)&1(*) &0.638(*) &0.587(*)\\
\hline
\hline
\multirow{4}{*}{VII} & 500    & {\bf 0.299} & 0.482(*) &0.343(*) &1.09(*) &0.818(*) &0.564(*) &0.583(*) \\
 & 1000    & {\bf 0.303} &0.465(*) &0.358(*) &1.51(*) &0.747(*) &0.589(*) &0.579(*) \\
 & 2000    & {\bf 0.309} &0.455(*) &0.383(*) &2(*) &0.966(*) &0.622(*) &0.578(*) \\
 & 3000     & {\bf 0.308} &0.46(*) &0.357(*) &2(*) &0.858(*) &0.626(*) &0.58(*)\\
\hline

\end{tabular}
\end{table}

\begin{table}[!h]
\caption{The average distance of the space estimated by each of the 7 methods tested to the true space $col(\vV)$ under various settings with $n=500$ and $p=6000$.
}\label{table:sim:3}
\begin{tabular}{|c|c|c|c|c|c|c|c|}
\hline
 & DT-SIR & SIRI-SIR & SIRS-SIR & SpSIR & SURE-SIR & DC-SIR & TP-SIR\\
\hline
\hline
I  & 0.694&0.631&{\bf 0.606}&2(*) &1.43(*) &0.97(*) &1.19(*) \\
\hline
II  &  0.446 &0.462&0.414&2(*) &1.74(*) &1.08(*) & {\bf 0.4} \\
\hline
III  & {\bf 1.35} & 1.56(*) &1.56(*) &2(*) &1.99(*) &1.88(*) &1.37 \\
\hline
IV   & 0.163 &{\bf 0.122}&0.245(*) &1.41(*) &0.27(*) &0.305(*) &0.195(*)\\
\hline
V  & 0.481(*) &{\bf 0.431}&0.486(*) &2(*) &1.62(*) &1.1(*) &0.995(*) \\
\hline
 VI   & 0.463(*) &{\bf 0.423} &0.494(*) &2(*) &1.62(*) &1.11(*) &0.999(*) \\
 \hline
VII  & 0.44 &{\bf 0.412} &0.477(*) &2(*) &1.61(*) &1.1(*) &1.03(*) \\
\hline
\end{tabular}
\end{table}

\begin{figure}[!h]
\centering
\includegraphics[height=60mm,width=50mm]{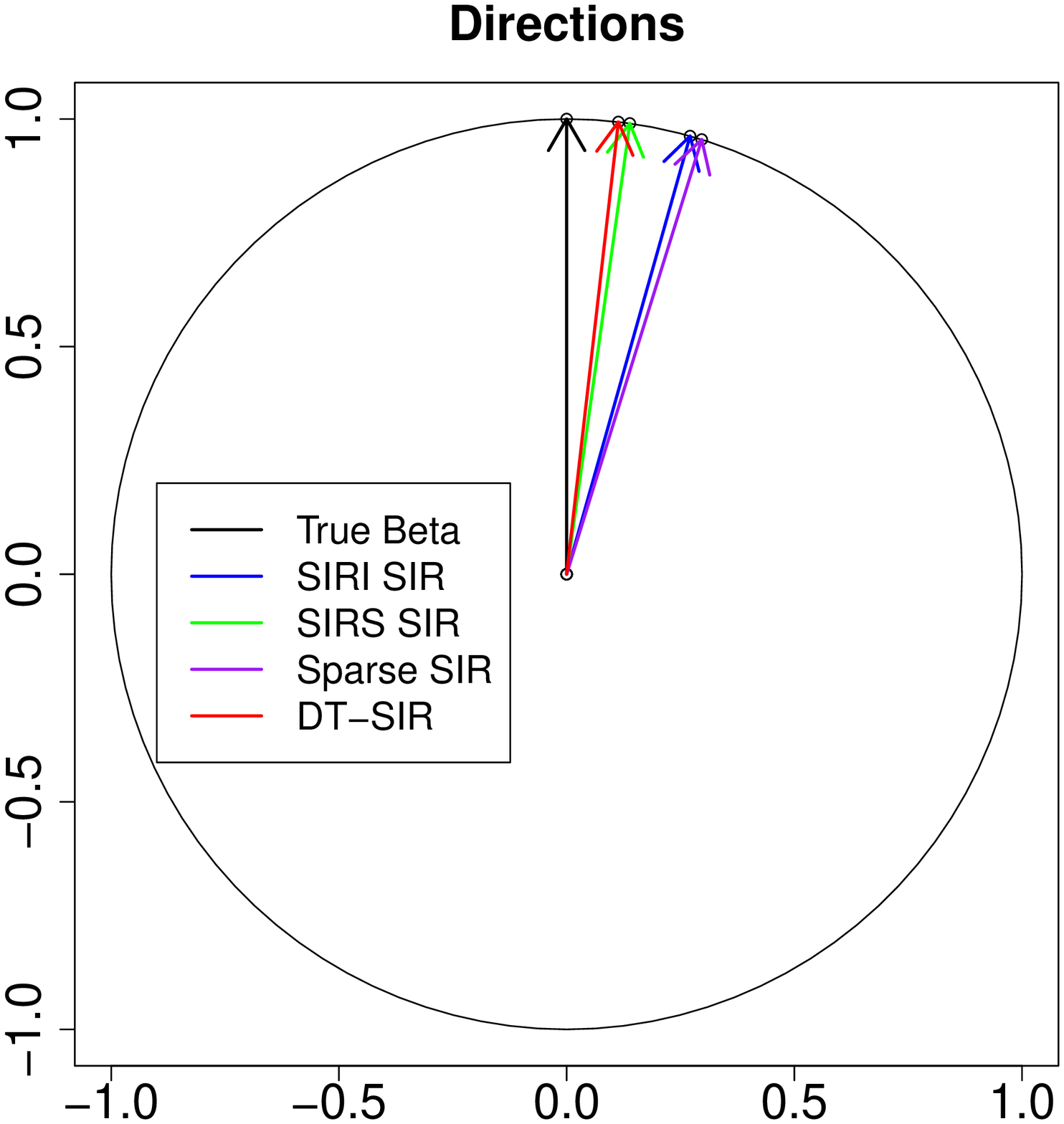}
\includegraphics[height=60mm,width=50mm]{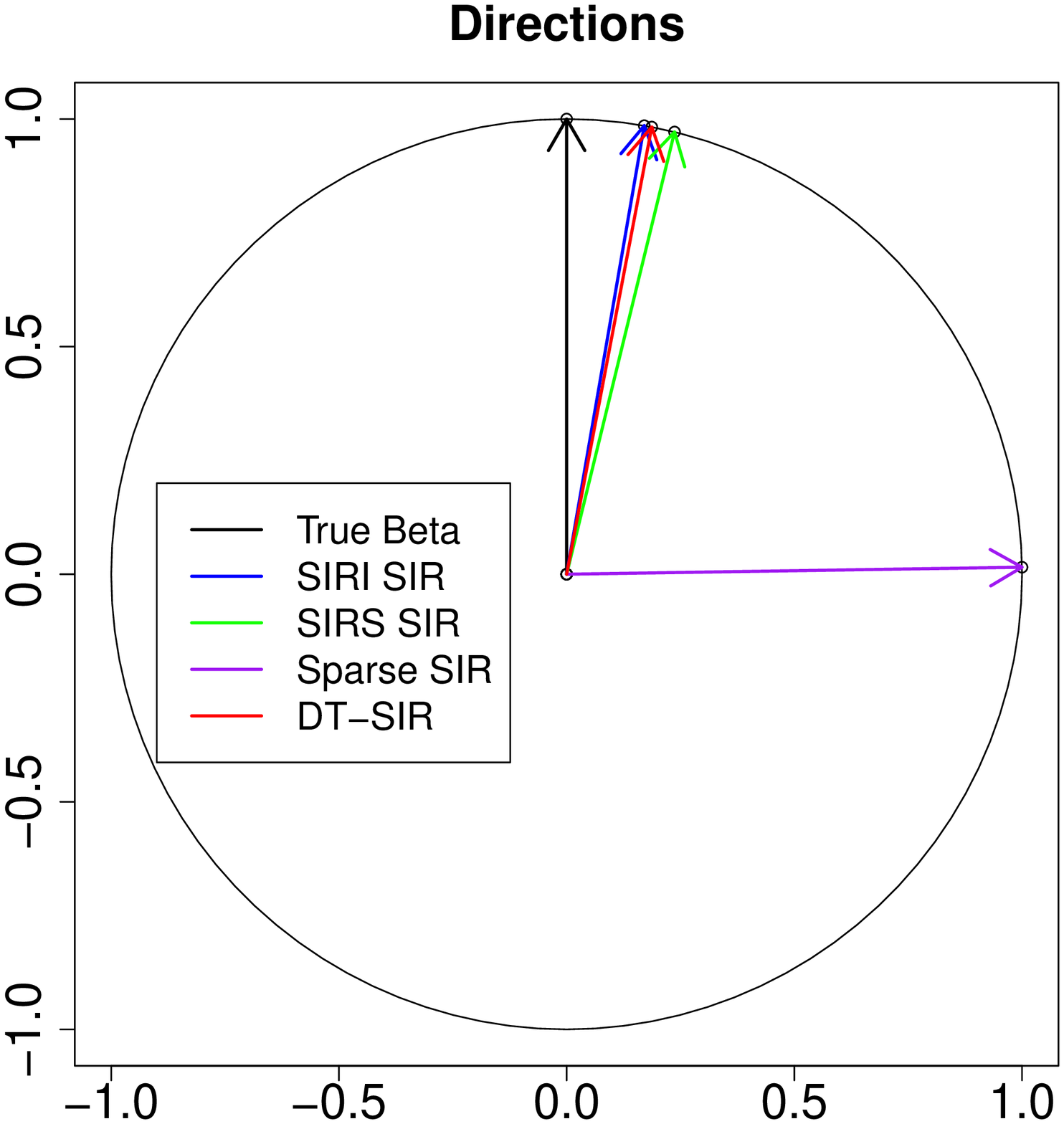}
\caption{ Simulated value of $E\angle(\hat{\bbeta}, \bbeta)$ for the various methods. Left panel: $(n,p)=(2000,1000)$; Right panel: $(n,p)=(500,1000)$.}\label{fig:2}
\end{figure}

To graphically show the performance of various methods, we consider setting IV with $d=1$.
Consider two cases when $(n,p)=(2000,1000)$ and $(n,p)=(500,100)$. We calculated the estimated directions $\hat{\bbeta}$ using various methods and computed the angle between $<\hat{\bbeta}>$ and $<\bbeta>$. 
We replicate this step 100 times to calculate the average angles for each method. The results are displayed in Figure \ref{fig:2}, which shows clearly that DT-SIR performed better than its competitors.

\begin{table}[!h]
\caption{Comparison of computing time under setting II.}\label{table:time}
\begin{tabular}{|c|c|c|c|c|c|c|c|}
\hline
& DT-SIR & SIRI-SIR & SIRS-SIR & SpSIR & SURE-SIR & DC-SIR & TP-SIR\\
\hline
\hline
n & \multicolumn{7}{c|}{p=1000} \\
\hline
500 & 1'' & 1'12'' & 7'' & 11'' & 1'' & 24''  & 29'' \\
\hline
1000& 2'' & 2'2'' & 20'' & 11'' & 1'' & 1'52''& 1'2''\\
\hline
2000& 3'' & 3'27'' & 1'14''& 13'' & 2'' & 7'38''& 2'18''\\
\hline
3000 & 4'' & 4'59'' & 2'45'' & 15'' & 3'' & 6'51''& 3'7''\\
\hline
\hline
p & \multicolumn{7}{c|}{n=2000} \\
\hline
500 & 1'' & 2'48'' & 35'' & 2'' & 1'' & 3'46''  & 1'7'' \\
\hline
1000& 3'' & 3'27'' & 1'14''& 13'' & 2'' & 7'38''& 2'18''\\
\hline
2000& 12'' & 4'55'' & 2'35'' & 1'39'' & 12'' & 14'24'' & 3'22''\\
\hline
3000 & 30'' & 6'0'' & 4'10'' & 5'19''& 30'' & 21'38'' & 6'17''\\
\hline
\end{tabular}
\end{table}

Additionally,  DT-SIR is computationally efficient. To show this, we report the computing time for one replication under Setting II for various pairs of $(n,p)$ in Table \ref{table:time}. All computations were done on a computer with Intel Xeon(R) E5-1620 CPU@3.70GHz and 16GB memory. It is clearly seen that DT-SIR performed as fast as SURE-SIR, and both were much faster than other competitors. Consider the case when $p=3000, n=2000$. The computation time of DT-SIR was only 30 seconds; while that for DC-SIR was 21 minutes and 38 seconds, and the that for TP-SIR was 6 minutes and 17 seconds.



\section{Conclusion.}\label{sec:conclusion}
When the dimension $p$ diverges to infinity, classical statistical procedures often fail unless additional structures such as sparsity conditions are imposed. 
Understanding boundary conditions of a statistical procedure provides us theoretical justification and practical guidance for our modeling efforts.
In this article, we provide a new framework to show that $\rho=\lim\frac{p}{n}$ is the phase transition parameter for the SIR procedure.
Under certain conditions, it is shown that the SIR estimator is consistent if and only if $\rho=0$.
When $\rho>0$, where the original SIR fails to be consistent, we  propose a two-stage method, DT-SIR for variable screening and selection in ultra-high dimension situations and show that the method is consistent.
We have used simulated examples to demonstrate the advantages of DT-SIR compared to its competitors.
This method is computationally fast and can be easily implemented for large data sets.  



\setcounter{section}{0}
\setcounter{subsection}{0}
\setcounter{figure}{0}
\setcounter{table}{0}
\setcounter{page}{1}
\renewcommand\thesection{\Alph{section}}
\renewcommand\thefigure{\Alph{section}\arabic{figure}}
\renewcommand\thetable{\Alph{section}\arabic{table}}
\numberwithin{equation}{section}

\begin{appendices}

\vspace*{4mm}
\begin{center}
{\bf Appendices}
\end{center}
\vspace*{4mm}

In the following two sections we offer some details about our theoretical derivations, but some more tedious intermediate steps (organized as Lemmas 6-21) are deferred to the Supplemental Document to this article, which is available on line.

\section{The Key Lemma.}
The following lemma plays an important role in developing the high dimensional theory for sliced inverse regression. 
The proof of this key lemma is lengthy and technical. 
It will be helpful to keep in mind that $H$ and $\nu$ (if they are not constants) grow at very slow rate compared with $c$ and $n$ (e.g., polynomial of $\log(n)$).
Let 
$
\vm(y) =\bbE[\vx|y],
$
and 
$
\vx =\vm(y)+\epsilon.
$
Notations
$
\vm_{h,j}, ~\overline{\vm}_{h,\cdot},~ \overline{\overline{\vm}},\quad$ and  
$\bepsilon_{h,j},~ \overline{\bepsilon}_{h,\cdot},~ \overline{\overline{\bepsilon}}
$
are similarly defined as  $\vx_{h,j}$,   $\overline{\vx}_{h,\cdot}$ and $\vxbarbar$ that were introduced before.

\begin{lem}\label{lem:new_statistic:deviation} Let $\vx \in \mathbb{R}^{p}$ be a sub-Gaussian random variable which is upper exponentially bounded by $K$ (see Definition \ref{def:subgaussion:upper:norm}). For any unit vector $\bbeta \in \mathbb{R}^{p}$, let $\vx(\bbeta)=\langle\vx,\bbeta\rangle$ and $\vm(\bbeta)=\langle\vm,\bbeta\rangle=\bbE[\vx(\bbeta)\mid y]$, we have the following:
\begin{itemize}
\vspace*{3mm}
\item[i)] If $var(\vm(\bbeta))= 0$ , there exists positive constants $C_{1}, C_{2}$ and $C_{3}$  such that  for any $b =O(1)$ and sufficiently large $H$,  we have
\begin{equation*} \label{eqn:deviation:outSets}
\begin{aligned}
\bbP(var_{H}(\vx(\bbeta))>b) 
\leq  C_{1}\exp\left(-C_{2}\frac{nb}{H^{2}} +C_{3}\log(H) \right).
\end{aligned}
\end{equation*}

\vspace*{3mm}
\item[ii)] If $var(\vm(\bbeta))\neq 0$ , there exists positive constants $C_{1}, C_{2}$ and $C_{3}$ such that,  for any $\nu>1$,  we have
\[
|var_{H}(\vx(\bbeta))-var(\vm(\bbeta))|\geq \frac{1}{2\nu}var(\vm(\bbeta))
\]
with probability at most 
\[
C_{1}\exp\left( -C_{2}\frac{n~var(\vm(\bbeta))}{H^{2}\nu^{2}}+C_{3}\log(H)\right).
\]
where we choose $H$ such that $H^{\vartheta}>C_{4}\nu$ for some sufficiently large constant $C_{4}$.
\end{itemize}
\end{lem}

\subsection{Proof of Lemma \ref{lem:new_statistic:deviation} i) }\label{sec:sec}

If $\vm(\bbeta)= 0$ (or equivalently $var(\vm(\bbeta))= 0$), since
\begin{align*}
\bar{\bepsilon}_{h,\cdot}(\bbeta)^{2}=&\left(\frac{c-1}{c}\frac{1}{c-1}\sum_{i=1}^{c-1}\bepsilon_{h,i}(\bbeta)+\frac{1}{c}\bepsilon_{h,c}(\bbeta)\right)^{2} \\
\leq & 2 \left(\frac{1}{c-1}\sum_{i=1}^{c-1}\bepsilon_{h,i}(\bbeta)\right)^{2}+2\left(\frac{1}{c}\bepsilon_{h,c}(\bbeta)\right)^{2}
\end{align*}
for $h=1,...,H-1$ and $\bar{\bepsilon}_{H,\cdot}(\bbeta)=\frac{1}{c}\sum_{i=1}^{c}\bepsilon_{H,i}(\bbeta)$,
we have
\begin{align*}
&var_{H}(\vx(\bbeta))-var(\vm(\bbeta))\\
=&\frac{1}{H}\sum_{h}^{H-1}\bar{\bepsilon}_{h,\cdot}(\bbeta)^{2}+\frac{1}{H}\bar{\bepsilon}_{H,\cdot}(\bbeta)^{2}\\
\leq & 
 \frac{2}{H}\left(\sum_{h}^{H-1}\left(\frac{1}{c-1}\sum_{i=1}^{c-1}\bepsilon_{h,i}(\bbeta)\right)^{2}+\bar{\bepsilon}_{H,\cdot}(\bbeta)^{2}\right)+\frac{2}{Hc^{2}}\sum_{h}^{H-1} \bepsilon_{h,c}(\bbeta)^{2} \\
 \triangleq& 2I+2II.
\end{align*}
Thus
\begin{align}
\bbP(var_{H}(\vx(\bbeta))>b) 
&\leq \bbP(I>b/4)+\bbP(II>b/4).
\end{align}

 Lemma \ref{cor:key_deviation:lemma} (iii) in Supplement implies that
\[
\bbP\left( \bepsilon(\bbeta)|_{y \in S_{h}} > t \right) \leq CH\exp\left(-\frac{t^{2}}{K^{2}}\right)
\]
for some positive constant $C$.  
Since $\bbE[\vx(\bbeta)|y]=0$,we have $\bbE[\vx(\bbeta)|y\in S_{h}]=0$.
From Lemma \ref{rmk:key_rmk}, we know that for $1\leq h \leq H-1$, $\bepsilon_{h,i}(\bbeta)$ can be treated as  $c-1$ $i.i.d.$ samples from $\bepsilon(\bbeta)|_{y\in S_{h}}$.  According to Lemma \ref{cor:key_deviation:lemma} (iv),
\[
\bbP\left( |\frac{1}{c-1}\sum_{i=1}^{c-1}\bepsilon_{h,i}(\bbeta)|>\sqrt{b}/2  \right)
\leq C_{1} \exp \left( \frac{-b(c-1)}{8C_{2}HK^{2}+4\sqrt{b}K} \right).
\]
Similarly, we have
\[
\bbP\left( |\frac{1}{c}\sum_{i=1}^{c}\bepsilon_{H,i}(\bbeta)|>\sqrt{b}/2  \right)
\leq C_{1} \exp \left( \frac{-bc}{8C_{2}HK^{2}+4\sqrt{b}K} \right).
\]
Thus, if $b=O(1)$ and $H$ is sufficiently large,  we have
\begin{align*}
\bbP(I>\frac{b}{4})\leq& C_{1}\left((H-1) \exp \left( \frac{-b(c-1)}{8C_{2}HK^{2}+4\sqrt{b}K} \right) + \exp \left( \frac{-bc}{8C_{3}HK^{2}+4\sqrt{b}K} \right)\right) \\
\leq & C_{1}\exp\left(-C_{2}\frac{cb}{H} +C_{3}\log(H) \right)
\end{align*}
for some positive constants $C_{1}, C_{2}$ and $C_{3}$. 

Since $\bepsilon_{i}(\bbeta)$ are i.i.d. samples from  a sub-Gaussian distribution $\bepsilon(\bbeta)$ with mean 0 and upper-exponentially bounded by $2K$.  Lemma \ref{lem:Hanson_Wright} implies that if $b=O(1)$ and $H$ is sufficiently large, we have
\begin{align*}
\bbP(II>b/4) \leq& \bbP(\frac{1}{n}\sum_{i}\bepsilon_{i}(\bbeta)^{2}>bc/4)\\
\leq & \bbP(\frac{1}{n}\sum_{i}\bepsilon_{i}(\bbeta)^{2}-
\bbE[\bepsilon(\bbeta)^{2}]>bc/4-\bbE[\bepsilon(\bbeta)^{2})\\
\leq & \bbP\left( \Big| \frac{1}{n}\sum_{i}\bepsilon_{i}(\bbeta)^{2}-\bbE[\bepsilon(\bbeta)^{2}]\Big| \geq cb/4-4K^{2} \right)  \\
\leq& C_{1}\exp\left( -C_{2}\frac{\sqrt{n}(cb/4-4K^{2})}{K^{2}} \right) \\
\leq & C_{1}\exp\left(-C_{2}\frac{cb}{H}+C_{3}\log(H) \right)
\end{align*}
for some positive constants $C_{1}, C_{2}$ and $C_{3}$ if $H$ is sufficiently large.
We used in above the fact that $\bbE[\bepsilon(\bbeta)^{2}]\leq 4K^{2}$.

To summarize, if $b=O(1)$ and $H$ is sufficiently large,  we have
\[
\bbP(var_{H}(\vx(\bbeta))>b) 
\leq  C_{1}\exp\left(-C_{2}\frac{cb}{H} +C_{3}\log(H) \right)
\]
for some  positive absolute constants $C_{1}, C_{2}$ and $C_{3}$.

\subsection{Proof of Lemma \ref{lem:new_statistic:deviation} ii) }

Since $\vx$ is sub-Gaussian and $\bbeta$ is unit vector, we know that $var(\vm(\bbeta))= O(1)$. If  $\vm(\bbeta)\neq 0$ (or equivalently $var(\vm\left(\bbeta)\right)\neq 0$), we have
\begin{align*}
&\Big|var_{H}(\vx(\bbeta))-var(\vm(\bbeta))\Big|\\
 =&\Big|\frac{1}{H}\sum_{h}\overline{\vx}_{h,\cdot}(\bbeta)^{2}-var(\vm(\bbeta))\Big|\\
=&  \Big| \frac{1}{H}\sum_{h}\overline{\vm}_{h,\cdot}(\bbeta)^{2}+\frac{2}{H}\sum_{h}\overline{\vm}_{h,\cdot}(\bbeta)\overline{\bepsilon}_{h,\cdot}(\bbeta)+\frac{1}{H}\sum_{h}\overline{\bepsilon}_{h,\cdot}(\bbeta)^{2}\\
&-var(\vm(\bbeta)) \Big| \\
\leq &A_{1}+A_{2}+A_{3}+A_{4},
\end{align*} 
 
where
\begin{equation}\label{eqn:five_terms}
\begin{aligned}
A_{1}&=\Big|\frac{1}{H}\sum_{h}\mu_{h}(\bbeta)^{2}-var(\vm(\bbeta))\Big|,\\
A_{2}&=\frac{1}{H}\sum_{h}\Big|\overline{\vm}_{h,\cdot}(\bbeta)^{2}-\mu_{h}(\bbeta)^{2}\Big|, \\
A_{3}&=\frac{1}{H}\sum_{h}\bar{\bepsilon}_{h,\cdot}(\bbeta)^{2},\\
A_{4}&=(\frac{1}{H}\sum_{h}\overline{\vm}_{h,\cdot}(\bbeta)^{2})^{1/2}(\frac{1}{H}\sum_{h}\overline{\bepsilon}_{h,\cdot}(\bbeta)^{2})^{1/2}.
\end{aligned}
\end{equation}
Lemma \ref{lem:new_statistic:deviation} ii) is a direct corollary of the following properties of $A_{i}$'s.

\begin{lem}\label{lem:tail_detail} Let the $A_{i}$'s be defined as in equation \eqref{eqn:five_terms}.  There exist positive constants $C_{1}$, $C_{2}$ and $C_{3}$, such that for any $\nu>1$  and $H$ satisfying $H^{\vartheta}=N_{1}\nu$ for sufficiently large $N_{1}$, we have that  each of the following events
\begin{itemize}
\item [{\bf i)}] 
$\Theta_{1}=\Big\{~ A_{1}\leq \frac{1}{4\nu} var(\vm(\bbeta)) ~\Big\},$

\vspace*{3mm}
\item [{\bf ii)}]  
$\Theta_{2}=\Big\{~A_{2} \leq \frac{1}{8\nu} var(\vm(\bbeta)) ~\Big\},$

\vspace*{3mm}
\item [{\bf iii)}] 
$
\Theta_{3}=\Big\{~A_{3}\leq \frac{1}{16\nu} var(\vm(\bbeta)) ~\Big\},
$

\vspace*{3mm}
\item [{\bf iv)}] 
$
\Theta_{4}=\Big\{~A_{4}  \leq \frac{1}{16\nu} var(\vm(\bbeta)) ~\Big\},
$
\end{itemize}
occurs with probability  at least 
\begin{equation}\label{tempppp1}
1-C_{1}\exp\left( -C_{2}\frac{c~var(\vm(\bbeta))}{H\nu^{2}}+C_{3}\log(H)\right).
\end{equation}
\end{lem}
\epf

\subsubsection{Proof of Lemma \ref{lem:tail_detail}. }

\paragraph{Proof of  ~${\bf i)} : $}
Recall definitions of the random intervals $S_h, h=1,2,\cdots, H$ and random variable $\delta_h=\delta_h(\omega)=\int_{y\in S_h(\omega)}f(y)dy$. We have

 \begin{align*}
&\Big|\frac{1}{H}\sum_{h}\left(\bold{\mu}_{h}(\bbeta)\right)^{2}-var(\vm(\bbeta)) \Big| \\
\leq&\Big|var(\vm(\bbeta))-\sum_{h}\delta_h\left(\bold{\mu}_{h}(\bbeta)\right)^{2}\Big|+\Big|\frac{1}{H}\sum_{h}\left(\bold{\mu}_{h}(\bbeta)\right)^{2}  -\sum_{h}\delta_h\left(\mu_h(\bbeta)\right)^{2}\Big|
\\ \triangleq & B_{1}+B_{2}
\end{align*}

Let $\epsilon=\frac{1}{Hn_{0}+1}$ where $n_{0}=N_{2}\nu$ for some sufficiently large constant $N_{2}$ and let event $E(\epsilon)$ be defined as in Lemma \ref{lem:elementary:deviation} in Section \ref{sec:assist:lemma}, i.e.,
 $E(\epsilon)=\Big\{~\omega ~\Big|~ |\delta_{h}-\frac{1}{H}| > \epsilon, \forall h~\Big\}$.
For any $\omega \in E(\epsilon)^{c}$, we have
\begin{align}
\nonumber B_{1}
 = & \sum_{h}\delta_{h}(\omega)var(\vm(\bbeta)|y\in S_{h}(\omega))\\
 \leq & (\frac{1}{H}+\epsilon) \sum_{h}var(\vm(\bbeta)|y\in S_{h}(\omega))\label{eqn:inline:arbit3:bound} \\
\leq & (1+H\epsilon)\frac{\bgamma_{3}}{H^{\vartheta}}var(\vm(\bbeta))\label{eqn:inline:arbit1:bound} \\
\leq& \frac{2\bgamma_{3}}{N_{1}\nu}var(\vm(\bbeta)), \label{eqn:inline:arbit4:bound} 
 \end{align}
 where inequality \eqref{eqn:inline:arbit3:bound} follows from the fact that $\delta_h(\omega)\leq \frac{1}{H}+\epsilon$, inequality \eqref{eqn:inline:arbit1:bound} follows from the sliced stable condition \eqref{def:sliced_stable}   and  inequality \eqref{eqn:inline:arbit4:bound}
follows from the requirement that $H^{\vartheta}>N_{1}\nu$, and the fact
\begin{align}
B_{2}\nonumber \leq & \epsilon \sum_{h}\left(\bbeta^{\tau}\mu_h\right)^{2} = \sum_{h} \frac{\epsilon}{\delta_h}\delta_h\left(\bbeta^{\tau}\mu_h\right)^{2}\\
\leq &\frac{H\epsilon}{1-H\epsilon}\sum_{h}\delta_h\left(\bbeta^{\tau}\mu_{h}\right)^{2}
\label{eqn:inline:arbit2:bound} \\
\nonumber \leq & \frac{2}{N_{2}\nu}\sum_{h}\delta_h\left(\bbeta^{\tau}\mu_{h}\right)^{2}
 \end{align}
where inequality \eqref{eqn:inline:arbit2:bound} follows from the fact $\delta_{h}\geq \frac{1}{H}-\epsilon$.

From \eqref{eqn:inline:arbit1:bound}, we observe that
\begin{align}\label{eqn:temp:temp:temp:1}
\sum_{h}\delta_h\left(\bold{\mu}_{h}(\bbeta)\right)^{2}& \leq \left(1+\frac{2\bgamma_{3}}{N_{1}\nu}\right)var(\vm(\bbeta)).
\end{align}
Combining with \eqref{eqn:inline:arbit2:bound}, we then have
\begin{align*}
B_{2}\leq & \frac{2}{N_{2}\nu}\left(1+\frac{2\bgamma_{3}}{N_{1}\nu}\right)var(\vm(\bbeta)).
\end{align*}
So when $E(\epsilon)^{c}$ occurs, we have
\begin{align}
\nonumber &B_{1}+B_{2}
 \leq \left( \frac{2\bgamma_{3}}{N_{1}\nu}+\frac{2}{N_{2}\nu}\left(1+\frac{2\bgamma_{3}}{N_{1}\nu}\right)\right)var(\vm(\bbeta)). \label{eqn:temp:temp:2}
\end{align}
Note that $N_{1}$ and $N_{2}$ can be chosen sufficiently large so that 
\begin{align}
B_{1}+B_{2}\leq \frac{4\bgamma_{3}}{N_{1}\nu} var(\vm(\bbeta))\leq \frac{1}{4\nu} var(\vm(\bbeta)).
\end{align} 
Consequently, conditioning on $E(\epsilon)^{c}$ where $\epsilon=\frac{1}{HN_{2}\nu+1}$, if we choose $H^{\vartheta}>N_{1}\nu$, then 
\begin{equation}\label{eqn:donot_ask_me_1}
\Big|\frac{1}{H}\sum_{h}\left(\bold{\mu}_{h}(\bbeta)\right)^{2}-var(\vm(\bbeta)) \Big| \leq \frac{1}{4\nu}var(\vm(\bbeta)).
\end{equation}
Since $var(\vm(\bbeta)) =O(1)$, $H^{\vartheta}>N_{1}\nu$ and $\epsilon=\frac{1}{HN_{2}\nu+1}$, the desired probability bound follows from  Lemma \ref{lem:elementary:deviation}, i.e.,
\begin{align*}
\bbP(E(\epsilon))\leq &C_{1}\exp\left(-\frac{Hc+1}{32(Hn_{0}+1)^{2}}+\log(H^{2}\sqrt{Hc+1})\right)\\
\leq &
C_{1}\exp\left(-C_{2}\frac{c~var(\vm(\bbeta))}{H\nu^{2}}+C_{3}\log(H) \right).
\end{align*}
for some positive constants $C_{1}, C_{2}$ and $C_{3}$. 
\epf

\begin{rmk}\normalfont
From \eqref{eqn:donot_ask_me_1}, conditioning on $E(\epsilon)^{c}$, we obtain the following two inequalities 
\begin{align}
\frac{1}{H}\sum_{h}\left(\bold{\mu}_{h}(\bbeta)\right)^{2} &\leq \left(1+\frac{4\bgamma_{3}}{H^{\vartheta}} \right) var\left(\vm(\bbeta) \right) \label{eqn:sum_of_sqaure_of_slice_mean:1}
\end{align}
and 
\begin{align}
\frac{1}{H}\sum_{h}|\bold{\mu}_{h}(\bbeta)| \leq \left(\left(1+\frac{4\bgamma_{3}}{H^{\vartheta}} \right) var\left(\vm(\bbeta) \right) \right)^{1/2}. \label{eqn:sum_of_absolute_of_slice_mean:2}
\end{align}
In particular,  $\frac{1}{H}\sum_{h}\left(\bold{\mu}_{h}(\bbeta)\right)^{2}$ and  $\frac{1}{H}\sum_{h}|\bold{\mu}_{h}(\bbeta)|$ are bounded by $O_{P}(1)$.
\end{rmk}

\vspace*{4mm}
\paragraph{Proof of  ~${\bf ii)}:$} Denote $\frac{1}{c-1}\sum_{i=1}^{c-1}\vm_{h,i}(\bbeta)$ by $\overline{\vm}_{h}'(\bbeta)$ and $\overline{\vm}_{H,\cdot}(\bbeta)$ by $\overline{\vm}_{H}'(\bbeta)$,
we have 
\begin{align*}
A_{2}\nonumber \leq &\frac{1}{H}\sum_{h=1}^{H}\Big|\overline{\vm}'_{h}(\bbeta)^{2}-\mu_{h}(\bbeta)^{2}\Big|+\frac{1}{Hc^{2}}\sum_{h=1}^{H}\vm_{h,c}(\bbeta)^{2}\\
&+\frac{2(c-1)}{c} \left(\frac{1}{H}\sum_{h=1}^{H}\overline{\vm}'_{h}(\bbeta)^{2}\right)^{1/2}\left(\frac{1}{Hc^{2}}\sum_{h=1}^{H}\vm_{h,c}(\bbeta)^{2}\right)^{1/2}+
\frac{2}{Hc}\sum_{h=1}^{H}\mu_{h}(\bbeta)^{2}\\
\triangleq & ~I+II+III+IV
\end{align*}
Before we start proving this part, we need to introduce two events and bound their probabilities. First,
let
\begin{equation}
 E_{1}(N_{3},\nu)=\Big\{~ \eta(\bbeta)> \frac{1}{N_{3}\nu}\sqrt{var(\vm(\bbeta))}~\Big\}.
\end{equation}
where
$\eta(\bbeta)=\max_{1\leq h \leq H} \left\{  \Big|\overline{\vm}'_{h}(\bbeta)-\mu_{h}(\bbeta)\Big|\right\}. 
$
According to Lemma \ref{cor:key_deviation:lemma} (i), (iv) and Bonferroni's inequality, we have
\begin{align}
\bbP\left( E_{1}(N_{3},\nu) \right) \leq& 2H\exp\left( \frac{1}{(N_{3}\nu)^{2}}\frac{-(c-1)var(\vm(\bbeta))}{2CHK^{2}+\frac{2}{N_{3}\nu}\sqrt{var(m(\bbeta))}K}\right)\\
\leq& C_{1}\exp\left(-C_{2}\frac{c~var(\vm(\bbeta))}{H\nu^{2}} +C_{3}\log(H)\right) \label{temptemp}
\end{align}
for some positive constants $C_{1}$, $C_{2}$ and $C_{3}$. 
Second, let 
\[
E_{2}(N_{4},\nu)\triangleq\Big\{ ~ II > \frac{1}{N_{4}\nu}var(\vm(\bbeta))  ~\Big \},
\]
then
\begin{align*}
\bbP\left(E(N_{4},\nu)   \right) &\leq \bbP\left(\frac{1}{nc}\sum_{i}m^{2}_{i}> \frac{var(\vm(\bbeta))}{N_{4}\nu}\right) \\
\leq& C_{1}\exp\left( -C_{2}\sqrt{n}(c~\frac{var(\vm(\bbeta))}{\nu})-K^{2})\right)\\
\leq & C_{1}\exp\left(-C_{2}\frac{c~var(\vm(\bbeta))}{H\nu} +C_{3}\log(H)\right)
\end{align*}
for some positive constant $C_{1}$, $C_{2}$ and $C_{3}$. It is easily to see $E(N_{4},\nu) \subset E(N_{4},\nu^{2})$.

\vspace*{3mm}
\subparagraph{For I.} 
Conditioning on the event $E(\epsilon)^{c}\cap E_{1}(N_{3},\nu)^{c}$, combining with \eqref{eqn:sum_of_absolute_of_slice_mean:2}, we have
\begin{align*}
I \leq &\frac{1}{H}\sum_{h}\eta(\bbeta)(\eta(\bbeta)+2|\mu_{h}(\bbeta)|)\leq \eta(\bbeta)^{2}+\frac{2\eta(\bbeta)}{H}\sum_{h}|\mu_{h}(\bbeta)|\\
\leq & \left( \left(\frac{1}{N_{3}\nu}\right)^{2}+\frac{2}{N_{3}\nu}\left(1+\frac{4\bgamma_{3}}{H^{\vartheta}} \right)^{1/2} \right)var(\vm(\bbeta))\\
\leq& \frac{1}{32\nu}var(\vm(\bbeta))
\end{align*}
 if $N_{3}$ is sufficiently large .

\begin{rmk} \normalfont
From above, conditioning on the event $E(\epsilon)^{c}\cap E_{1}(N_{3},\nu)^{c}$, we have
\begin{align}
\frac{1}{H}\sum_{h=1}^{H}\overline{\vm}'(\bbeta)^{2} \leq& \frac{1+32\nu}{32\nu}var(\vm(\bbeta)).\label{eqn:number:87}
\end{align}
\end{rmk}


\vspace*{3mm}
\subparagraph{For II.} Conditioning on $E_{2}(N_{4},\nu)^{c}$, we have
$
II\leq \frac{var(\vm(y))}{N_{4}\nu}.
$

\vspace*{3mm}
\subparagraph{For III.} 
When the event $E(\epsilon)^{c}\cap E_{1}(N_{3},\nu)^{c}\cap E_{2}(N_{4},\nu^{2})^c$ occurs, according to equation (\ref{eqn:number:87}),
\begin{align}
\nonumber III 
 &\leq \frac{2(c-1)}{c}\sqrt{\frac{1+32\nu}{32\nu}}\frac{1}{\sqrt{N_{4}}\nu}var(\vm(\bbeta)) <\frac{1}{16\nu}var(\vm(\bbeta)) .
\end{align}
if $N_{4}$ is sufficiently large.

\vspace*{3mm}
\subparagraph{For VI.} When the event $E(\epsilon)^{c}\cap E_{1}\left(N_{3},\nu \right)^{c}\cap E_{2}(N_{4},\nu)^{c}$ occurs, from \eqref{eqn:donot_ask_me_1}, we know
\[
VI=\frac{2}{Hc}\sum_{h}\mu_{h}(\bbeta)^{2} \leq \frac{9}{4c}var(\vm(\bbeta)) <\frac{1}{16\nu}var(\vm(\bbeta)).
\]

\vspace{3mm} To summarize, we know that there exist positive constant $C_{1},C_{2}$, $C_{3}$ and $C_{4}$ such that 
\[
A_{2}\leq I+II+III+VI  \leq \frac{1}{8\nu} var(\vm(\bbeta))
\]
holds on the event   $E(\epsilon)^{c}\cap E_{1}(N_{3},\nu)^{c}\cap E_{2}(N_{4},\nu^{2})^{c}$ which is  with probability at least
\begin{align}
\nonumber 1&-C_{1}\exp\left( -C_{2}\frac{c~var(\vm(\bbeta))}{H\nu^{2}}+C_{3}\log(H)\right)
\end{align}
for some positive constants $C_{1}$,$C_{2}$ and $C_{3}$.

\paragraph{Proof of  ~${\bf iii)} : $}
Similar to the proof of Lemma \ref{lem:new_statistic:deviation} (i) we have
\[
\bbP(A_{3}>b) \leq C_{1}H \exp \left( \frac{-(c-1)b}{8C_{2}HK^{2}_{1}+4\sqrt{b}K_{2}} \right)
\]
for some positive constants $C_{1}$, $C_{2}$ and $C_{3}$.
In particular , if we take $b=\frac{1}{16\nu} var(\vm(\bbeta))$, we know that
\[
A_{3}\leq \frac{1}{16\nu}var(\vm(\bbeta))
\]
with probability at least
\[
1-C_{1}\exp\left(-C_{2}\frac{c~var(\vm(\bbeta))}{H\nu^{2}} +C_{3}\log(H)\right)
\]
for some positive constant $C_{1}$, $C_{2}$ and $C_{3}$ .

\paragraph{Proof of  ~${\bf iv)} : $} Let 
\begin{align}
\nonumber D_{1}\triangleq\frac{1}{H}\sum_{h}\overline{\vm}_{h,\cdot}(\bbeta)^{2},\quad 
D_{2}\triangleq A_{3}=\frac{1}{H}\sum_{h}\bar{\bepsilon}_{h,\cdot}(\bbeta)^{2}
\end{align}
Consequently,
\begin{align}
\nonumber  &\bbP\left( D^{1/2}_{1}D^{1/2}_{2} > \frac{1}{16\nu}var(\vm(\bbeta))\right)\\
\leq & \bbP\left( |D_{1}|>\frac{2\nu+1}{2\nu}var(\vm(\bbeta))\right)+\bbP\left( D_{2} > \frac{var(\vm(\bbeta))}{(2\nu+1)16\nu} \right)\label{eqn:b1:b2}
\end{align}

Note that
\begin{align}
\nonumber | D_{1}-var(\vm(\bbeta))|
\leq A_{2}+A_{1}
\end{align}
According to (i) and (ii),  the right hand side of (\ref{eqn:b1:b2}) is bounded by
\begin{align*}
C_{1}\exp\left(-C_{2}\frac{c~var(\vm(\bbeta)}{H\nu^{2}}+C_{3}\log(H) \right)
\end{align*}
for some positive constants $C_{1},C_{2}$ and $C_{3}$.

\epf

\section{Proofs of theorems in section \ref{sec:sir}.} \label{sec:appendix:A}

\subsection{Proof of Theorem \ref{thm:high:dim}.}


Let $\mathcal{S}$ be the central subspace of dimension $d \ll p$,
i.e., $y\independent \vx | \vP_\mathcal{S}\vx$ and $dim(\mathcal{S})=d$.  We have the decomposition 
\begin{equation}\label{eqn:decomposition}
\begin{aligned}
\vx &= \vP_{\mathcal{S}}\vx + \vP_{\mathcal{S^\perp}}\vx \triangleq \vz + \vw \\
&=\bbE[\vz|y]+\vz-\bbE[\vz|y]+\vw\triangleq\vm+\vv+\vw
\end{aligned}
\end{equation}
where $\vz=\vP_{\mathcal{S}}\vx$, $\vm=\bbE[\vz|y]$, $\vv=\vz-\bbE[\vz|y]$ and $\vw= \vP_{\mathcal{S^\perp}}\vx$.  Note that $\vm$ lies in the central curve, $\vv$ lies in the central space and $\vw$ lies in the space perpendicular to $\mathcal{S}$.
We introduce 
\begin{align}
\vm_{h,j}, ~\overline{\vm}_{h,\cdot},~ \overline{\overline{\vm}},\quad 
\vz_{h,j},~ \overline{\vz}_{h,\cdot},~ \overline{\overline{\vz}},\mbox{ and }
\vw_{h,j},~ \overline{\vw}_{h,\cdot},~ \overline{\overline{\vw}}
\end{align}
similar to the definition of $\vx_{h,j}$,   $\overline{\vx}_{h,\cdot}$ and $\vxbarbar$.  
Consequently, we can define $\widehat{\bLambda}_{\vz}$ and have the following decomposition
\begin{equation}\label{eqn:decomposition:model}
\widehat{\bLambda}_{H}\equiv \frac{1}{H}\sum_{h}\overline{\vx}_{h,\cdot}\overline{\vx}_{h,\cdot}^{\tau}
= \widehat{\Lambda}_{\vz}+\mathcal{Z}\mathcal{W}^{\tau}
+\mathcal{W}\mathcal{Z}^{\tau}+\mathcal{W}\mathcal{W}^{\tau},
\end{equation}
where 
\begin{align*}
\mathcal{Z}=\frac{1}{\sqrt{H}}\left(\overline{\vz}_{1,\cdot},...,\overline{\vz}_{H,\cdot}\right) \mbox{ and } \mathcal{W}=\frac{1}{\sqrt{H}}\left(\overline{\vw}_{1,\cdot},...,\overline{\vw}_{H,\cdot}\right).
\end{align*}

We need to bound $\|\widehat{\bLambda}_{\vz}-\bLambda_{\vp}\|_{2} $  and $\|\mathcal{W}\mathcal{W}^{\tau}\|_{2}$.

\begin{lem}\label{lem:pure:error}
\begin{align}
\|\mathcal{W}\mathcal{W}^{\tau}\|_{2}\leq O_{P}(\frac{H^{2}p}{n}) 
\end{align}
\proof
From Lemma \ref{lem:new_statistic:deviation},
for any unit vector $\bbeta \perp col(\bLambda)$, i.e. $var(\vm(\bbeta))=0$, we have
\begin{align}
\bbP(\bbeta^{\tau}\mathcal{W}\mathcal{W}^{\tau}\bbeta> C\frac{H^{2}p}{n}) \leq C_{1}\exp\left(-C_{2}p+\log(H) \right).
\end{align}
for some positive constants $C_{1}$ and $C_{2}$.
Then the $\varepsilon$-net argument (see e.g., \cite{Vershynin:2010}) implies that $\|\mathcal{W}\mathcal{W}^{\tau}\|\leq O_{P}(\frac{H^{2}p}{n})$  \epf
\end{lem}

\begin{lem}\label{lem:eigenvalue}
\begin{align} 
\|\widehat{\bLambda}_{\vz}-\bLambda_{p}\| \leq O_{P}\left(\frac{1}{H^{\vartheta}}\right). 
\end{align}
As a direct corollary, we have $\|\widehat{\bLambda}_{\vz}\|\leq O_{P}(1)$.
\proof From Lemma \ref{lem:new_statistic:deviation}, we have
\begin{align}
\bbP\left( \big|\bbeta^{\tau}(\widehat{\bLambda}_{\vz}-\bLambda) \bbeta\big|> \frac{C}{H^{\vartheta}}\|\bLambda\|_{2} \right)
\leq C_{1}\exp\left( -C_{2}\frac{c~var(\vm(\bbeta))}{H^{1+2\vartheta}}+C_{3}\log(H) \right).\nonumber
\end{align}
Note that we only need to verify it for $\bbeta \in col(\bLambda_{p})$, which is a $d$-dimensional space. Then the $\varepsilon$-net argument implies that $\|\widehat{\bLambda}_{\vz}-\bLambda_{p}\|_{2} \leq O_{P}\left(\frac{1}{H^{\vartheta}} \right)$.
\epf 
\end{lem}

Theorem \ref{thm:high:dim} follows from Lemma~\ref{lem:eigenvalue} and Lemma \ref{lem:pure:error}. In fact,
\begin{align*}
\|\widehat{\bLambda}_{H}-\bLambda_{p}\| &\leq \|\widehat{\bLambda}_{\vz}-\bLambda_{p}\|+\|\mathcal{Z}\mathcal{W}^{\tau}
+\mathcal{W}\mathcal{Z}^{\tau}\|_{2}+\|\mathcal{W}\mathcal{W}^{\tau}\|_{2}\\
& \leq O_{P}\left(\frac{1}{H^{\vartheta}}+\sqrt{\frac{H^{2}p}{n}} +\frac{H^{2}p}{n}\right).
\end{align*}
\epf

\subsection{Proof of Theorem \ref{thm:space:consistency}. } 

Theorem \ref{thm:space:consistency} is a direct corollary of Theorem \ref{thm:high:dim} and Lemma \ref{lem:folklore}. In fact, we have:
\begin{align*}
&\|\widehat{\bSigma}^{-1}_{\vX}\widehat{\bLambda}_{H}-\bSigma^{-1}_{X}\bLambda_{p}\|_{2}\\
&\leq \|\widehat{\bSigma}^{-1}_{\vX}-\bSigma^{-1}_{\vX}\|_{2}\|\widehat{\bLambda}_{H}\|_{2}+\|\bSigma^{-1}_{\vX}\|_{2}\|\widehat{\bLambda}_{H}-\bLambda_{p}\|_{2},
\end{align*}
which  $\rightarrow 0$ if $\rho=\lim_{n \rightarrow \infty}\frac{p}{n} = 0$.
\epf \\

\subsection{Proof of Theorem \ref{thm:main1.1}.} \ \\

{\bf (i)} \  
The proof for part (i) is similar to the proof of Theorem \ref{thm:high:dim} and the standard Gaussian assumption on $\vx$ simplifies the argument and improves the results. 
Since $\vw=\vP_{\mathcal{S}^{\perp}}\vx$ is normal and  independent of $y$, there exists a normal random variable $\bepsilon\sim N(\vzero, \vI)$ such that $\vw=\bSigma^{1/2}\bepsilon$ where $\bSigma=cov(\vw)$. Using the decomposition \eqref{eqn:decomposition:model},  we may write
\begin{align}
\mathcal{W}=\frac{1}{\sqrt{Hc}}\bSigma^{1/2}\vE_{p\times H}
\end{align}
where $\vE_{p,H}$ is a $p\times H$ matrix with $i.i.d.$ standard normal entries.
Corollary \ref{cor:final:done} implies that
\[
\|\mathcal{W}\mathcal{W}^{\tau}\|_{2}\leq C\left(\sqrt{\frac{p}{n}}+\sqrt{\frac{H}{n}} \right)^{2}\leq O_{P}\left(\frac{p}{n}\right).
\]
Lemma \ref{lem:eigenvalue} implies
\[
\|\widehat{\bLambda}_{\vz}\|_{2}\leq \|\bLambda_p\|_{2}+O_{P}\left(\frac{1}{H^{\vartheta}}\right).
\]
By the Cauchy inequality, we have
\begin{align*}
\|\mathcal{C}_{1}\|^{2}_{2}&\leq \|\widehat{\bLambda}_{\vz}\|_{2}\|\mathcal{W}\mathcal{W}^{\tau}\|_{2}
\leq O_{P}\left(\frac{p}{n} \right).
\end{align*}
Thus, 
\[
\|\widehat{\bLambda}_{H}-\bLambda_p\|_{2} \leq  O_{P}(\frac{1}{H^{\vartheta}}+\frac{p}{n}+\sqrt{\frac{p}{n}}).
\]
In particular, if $H, n\rightarrow \infty$ and $\rho=\lim\frac{p}{n}\in (0,\infty)$,  we know that  $\|\widehat{\bLambda}_{H}-\bLambda_{p}\|_{2}$ is  dominated by  $\rho\vee\sqrt{\rho}$ as a function of $\rho$. \epf \\

{\bf (ii)}  \ The proof for part (ii) is similar to the proof of Theorem 2 in \cite{Johnstone:Lu:2009} but is technically more challenging. 
Let  
$
D=\mathcal{Z}\mathcal{Z}^{\tau}+\mathcal{W}\mathcal{W}^{\tau} $ 
and 
$ 
B=\mathcal{Z}\mathcal{W}^{\tau}
+\mathcal{W}\mathcal{Z}^{\tau}, 
$
then
\begin{align*}
\widehat{\bLambda}_{H}
=&D+B.
\end{align*}
Since we are working on single index model with $\vx$ is standard normal,  $\vz=P_{\bbeta}\vx=\bbeta z(y)$ for some scalar function $z(y)$ and $\vw=P_{\bbeta^{\perp}}\vx$ are independent normal random variables.  
Let $\bSigma=var(\vw)$, then we can write
\[
\mathcal{W}=\frac{1}{\sqrt{Hc}}\bSigma^{1/2}\bold{E}
\]
where $\bold{E}$ is a $p\times H$ matrix with $i.i.d.$ standard normal entries.

Since $\vz=\bbeta z(y)$, we have
$\mathcal{Z}=\frac{1}{\sqrt{H}}\bbeta(\overline{z}_{1,\cdot},\overline{z}_{2,\cdot},...,\overline{z}_{H,\cdot})$. To ease notation, let $\btheta^{\tau}=(\overline{z}_{1,\cdot},\overline{z}_{2,\cdot},...,\overline{z}_{H,\cdot})$, then

\begin{equation}\label{eqn:middle_notation}
\begin{aligned}
D&=\frac{1}{H}\|\btheta\|^{2}\bbeta\bbeta^{\tau}+\frac{1}{n}\bSigma^{1/2} \vE\vE^{\tau}\bSigma^{1/2}\\
 B&=\bbeta \vu^{\tau}+\vu\bbeta^{\tau} \mbox{ where \quad }\vu=\frac{1}{H\sqrt{c}}\bSigma^{1/2} \vE \btheta   .
\end{aligned}
\end{equation}

Let $0<\alpha < arctan(\frac{1}{16})$ and 
\begin{align} \label{angle_neighbour}
N_{\alpha}=\Big\{~\vx \in \mathbb{R}^{p} :\angle\left(\vx,\bbeta\right) \leq \alpha \mbox{ and } \|x\|=1 ~\Big\}
\end{align}
be the set of unit vectors making angle at most $\alpha$ where $\angle(\vx,\vy)$ is the angle between the vectors $\vx$ and $\vy$.
In order to proceed, we need the following lemma.

\begin{lem}\label{angle_lemma}
Let $\widehat{\bbeta}$ and $\widehat{\bbeta}_{-}$ be the principal eigenvector of  $S_{+}\triangleq D+B$ and $S_{-}\triangleq D-B$, respectively. There exists a positive constant $\omega(\alpha)$ such that for any $\widehat{\bbeta}\in N_{\alpha}$, i.e., $\angle\left(\widehat{\bbeta},\bbeta\right) \leq \alpha$, we have  
\begin{align}\label{high_probability}
\angle\left(\widehat{\bbeta},\widehat{\bbeta}_{-}\right)\geq \frac{1}{3}\omega(\alpha)
\end{align}
 with probability converging to one as $n\rightarrow \infty$.
 \proof The proof is presented in \cite{lin2015sir}.
\end{lem}

Note that $S_{+}$ and $S_{-}$ have the same distribution (viewed as functions of random terms $\vE$ and $\theta$):
\[
S_{-}(\vE,\theta)=S_{+}(-\vE,\theta).
\]
Let $\mathcal{A}_{\alpha}$ denote the event $\Big\{\angle\left(\widehat{\bbeta},\bbeta\right) \leq \alpha\Big\}\cup \Big\{\angle\left(\widehat{\bbeta}_{-},\bbeta\right) \leq \alpha\Big\}$, then 
\begin{align*}
\bbE[\angle\left(\widehat{\bbeta},\bbeta\right)] & \geq \bbE[\angle\left(\widehat{\bbeta},\bbeta\right), \mathcal{A}^{c}_{\alpha}]  +\bbE[\angle\left(\widehat{\bbeta},\bbeta\right), \mathcal{A}_{\alpha}]  \\
&\geq \bbE[\angle\left(\widehat{\bbeta},\bbeta\right), \mathcal{A}^{c}_{\alpha}] +\frac{1}{2}\bbE[\angle\left(\widehat{\bbeta},\widehat{\bbeta}_{-}\right), \mathcal{A}_{\alpha}] \\
&\geq \min\{\alpha,\frac{\omega(\alpha)}{6}\} >0.
\end{align*}\epf








\section{Proof of Lemma \ref{angle_lemma}.} We need the following lemmas.

\begin{lem} \label{temp:upperbound}
Recall that $\vu=\frac{1}{H\sqrt{c}}\bSigma^{1/2} \vE \btheta $ defined as in \eqref{eqn:middle_notation}, then
there exist positive constants $C_{1}$ and  $C_{2}$ such that
\[
0<C_{1} \leq \|\vu\|_{2} \leq C_{2} 
\]
with probability converging to one as $n \rightarrow \infty$.
\end{lem}

\begin{lem}\label{lem:ass:geometry_I}
Assuming conditions in Theorem \ref{thm:main1.1}, let B  and 
$N_{\alpha}$ be defined as in \eqref{eqn:middle_notation} and\eqref{angle_neighbour} respectively where $0<\alpha < arctan(\frac{1}{16})$ .  
\begin{itemize}
\item[i)]  There exists positive constant $C_{1}$ such that for any $\vx \in N_{\alpha}$, we have $\|B\vx\|\geq C_{1}$  with probability converging to one as $n \rightarrow \infty$;
\item[ii)] For any $\vx \in N_{\alpha}$, we have $\Big|\cos\angle(\vx,B\vx) \Big| \leq 4 \alpha$  with probability converging to one as $n\rightarrow \infty$.
\end{itemize}
\end{lem}

The following lemma is borrowed from \cite{johnstone2004sparse}.
\begin{lem}\label{lem:ass:geometry_II}
Let $\bxi$ be a principal eigenvector of a non-zero symmetric  matrix M. For any $\boldeta\neq 0$,
\[
\angle(\boldeta, M\boldeta)\leq 3 \angle(\boldeta, \bxi).
\]
\end{lem}

The proof of Lemma \ref{angle_lemma} is made plausible by reference to the Figure \ref{fig:assistant}.

\begin{figure}[H]
\includegraphics[natwidth=50mm]{./triangle.png}
\caption{An illustrated graph}\label{fig:assistant}
\end{figure}

\vspace{3mm}
Since
\begin{align}\label{lowerbound_angle}
\sin\left(\angle\left(\widehat{\bbeta},S_{-}\widehat{\bbeta}  \right)\right)&=\sin\left(\omega_{1}+\omega_{2}\right)=\sin\left(  \pi -\omega_{1}-\omega_{2}\right)\\
&\geq \min\Big\{~\sin\left(\omega_{1}\right), \sin\left(\omega_{3}\right)~\Big\},
\end{align}
we only need to prove that there exists a positive small constant $\omega(\alpha)$ ( $< \frac{\pi}{2}$ ) such that $\sin\left(\omega_{1}\right), \sin\left(\omega_{2}\right)\geq \sin\left(\omega(\alpha)\right) $ . In fact, if such $\omega(\alpha)$ exists, we may choose $\vM=S_{-}, \xi=\widehat{\bbeta}_{-}$ in  Lemma \ref{lem:ass:geometry_II} and get 
\[
\angle\left(\widehat{\bbeta},\widehat{\bbeta}_{-}\right)\geq \frac{1}{3}\angle\left(\widehat{\bbeta},S_{-}\widehat{\bbeta}\right) \geq \frac{1}{3}\omega(\alpha). 
\]

\subparagraph{For $\omega_{3}$.}  From Lemma \ref{lem:ass:geometry_I}  ii), $\Big|\cos\angle(\widehat{\bbeta},B\widehat{\bbeta})\Big|\leq 4\alpha$, we know that there exists positive constants $\delta(\alpha)(<\frac{\pi}{2})$ such that $\sin\omega_{3}\geq sin\left(\delta(\alpha)\right)$ . 

\vspace*{3mm}
\subparagraph{For $\omega_{1}$.}  Applying the law of sines to the triangle $\triangle(O,O_{1},O_{2})$ , we have
\begin{align}\label{low_of_since}
\frac{\sin\omega_{1}}{\|B\widehat{\bbeta}\|}=\frac{\sin\omega_{3}}{\|D\widehat{\bbeta}\|} \left( =\frac{\sin\angle\left(B\widehat{\bbeta},\widehat{\bbeta}\right)}{\|D\widehat{\bbeta}\|} \right).
\end{align}

Note that  from Lemma \ref{lem:ass:geometry_I} i),  there exists a constant $C_{1}>0$ such that $\|B\widehat{\bbeta}\|>C_{1}$ and 
 \[
 \|D\widehat{\bbeta}\|\leq \|D\| \leq \frac{1}{H}\|\theta\|^{2}+\|\frac{1}{n}EE^{\tau}\|,
 \]
is bounded by an absolute constant C given $\lim_{n\rightarrow \infty} \frac{p}{n}=\rho\neq 0$ and sliced stable condition.  Then \eqref{low_of_since} implies
\[
\sin \omega_{1}=\frac{\|B\widehat{\bbeta}\|\sin\angle(B\widehat{\bbeta},\widehat{\bbeta})}{\|D\widehat{\bbeta}\|} \geq \frac{C_{1}\sin\delta(\alpha)}{C} \geq \sin \omega' >0
\] 
where $\omega'$ ( $<\frac{\pi}{2}$ )is an small angle such that the last inequality holds. 
In particular, we have $\omega_{1}\geq \omega'$. Hence 
\begin{align*}
\angle(\widehat{\bbeta},S_{-}\widehat{\bbeta}) 
\geq & \omega'\wedge\delta(\alpha)\triangleq\omega(\alpha)
\end{align*}

\epf

\subsection{Proof of Lemma \ref{temp:upperbound}}
\textit{Proof} : In fact, let $\vT$ be an orthogonal matrix such that $\vT\bbeta=(1,0\cdots,0)^{\tau}$ and $\vM=\vT\bbeta\bbeta^{\tau}\vT^{\tau}$, then
\begin{align*}
cH^{2}\vu^{\tau}\vu
&=\btheta^{\tau}\vE^{\tau}\bSigma E\btheta\\
&=\btheta^{\tau}\vE^{\tau}\vE\btheta-\btheta^{\tau}\vE^{\tau}\bbeta\bbeta^{\tau}\vE\btheta\\
&=\btheta^{\tau}\vE^{\tau}\vT^{\tau}\vT\vE\btheta-\btheta^{\tau}\vE^{\tau}\vT^{\tau}(\vT\bbeta)\bbeta^{\tau}\vT^{\tau}\vT\vE\btheta\\
&\overset{d.}{=}\btheta^{\tau}\vE^{\tau}\vE\btheta-\btheta^{\tau}\vE^{\tau}\vM\vE\btheta\\
&\overset{d.}{=}\frac{p-1}{p}\btheta^{\tau}\vE^{\tau}\vE\btheta,
\end{align*}
where $\overset{d.}{=}$ means equal in distribution. Note that $\vE^{\tau}\vE$ is full rank $H\times H$ matrix,  combining with Lemma \ref{random:nonasymptotic} , we know that
\[
C_{1}(1-\sqrt{\frac{H}{p}})^{2} \leq\lambda_{min}(\frac{1}{p}\vE^{\tau}\vE)\leq \lambda_{max}(\frac{1}{p}\vE^{\tau}\vE) \leq C_{2}(1+\sqrt{\frac{H}{p}})^{2}
\]
for some positive constants $C_{1}$ and $C_{2}$ with probability at least $1-2\exp(-p/8)$.  
Note that $\lim\frac{p}{n}=\rho>0$ as $n \rightarrow \infty$ and $n=Hc$, we know there exists positive constants $C_{1}$ and $C_{2}$ such that
\[
C_{1} \frac{1}{H} \|\btheta\|^{2} \leq \|\vu\|^{2} \leq C_{2}\frac{1}{H} \|\btheta\|^{2}
\]
with probability at least $1-2\exp(-p/8)$.

On the other hand, the sliced stable condition  implies that $\lim\frac{1}{H}\|\btheta\|^{2}$ exists ($\not =0$), so  $\|\vu\|^{2}$ is bounded away from 0 and $\infty$ with probability 1 as $n\rightarrow \infty$. \epf

\subsection{Proof of Lemma \ref{lem:ass:geometry_I} }

\textit{Proof}: For $i)$,  $\forall \vx \in N_{\alpha}$, let
\[
\vx=\cos(\delta )\bbeta+\sin(\delta)\boldeta \mbox{ where } \boldeta \perp \bbeta, \|\boldeta\|=1, \delta \leq \alpha.
\]

Since $B\vx=\cos(\delta)\vu+(\vu^{\tau}\boldeta)\sin(\delta)\bbeta$, we have:
\begin{align*}
\|B\vx\|&\geq \cos(\delta) \|\vu\|-\sin(\delta) \|\vu\|\geq \frac{1}{2}\cos(\delta) \|\vu\|\\
&\geq \frac{1}{2}\cos(\alpha)\|\vu\|>\frac{C_{1}}{4}>0
\end{align*}
for some positive constant $C_{1}$.

\vspace{3mm}
For $ii)$, since
\[
\vx^{\tau}B\vx=2(\vu^{\tau}\boldeta)\cos(\delta)\sin(\delta) 
\]
we have that uniformly over $N_{\delta}$,
\[
|\vx^{\tau}B\vx|\leq |\vu^{\tau}\boldeta|\sin(2\delta)
\]
which in turn implies:
\[
\Big|\cos(\angle(B\vx,\vx))\Big|=\frac{|\vx^{\tau}B\vx|}{\|\vx\|\|B\vx\|}
\leq \frac{\sin(2\delta)|\vu^{\tau}\boldeta|}{\frac{1}{2}\cos(\delta)\|\vu\|}\leq 4 \delta \leq 4 \alpha.
\]
\epf

\section{Appendix B: Proofs of  Theorems \ref{thm:screen:consistent} to \ref{thm:operator_norem:high_dimensional}.} \label{sec:appendixC}

\subsection{Proof of Theorem \ref{thm:screen:consistent}.}
Let $\vx(k)=<\vx,\bbeta_{k}>$ where $\bbeta_{k}=(0,..,1,...0) \in \mathbb{R}^{p}$ with the only 1 at the $k$-th position.
Recall that  
 \begin{align*}
  \mathcal{T}=&\Big\{~k ~\Big|~ \bbE[\vx(k)|y] \mbox{ is not constant. }\Big\} \\
 \mathcal{I}_{p}(t)=&\Big\{~k~\Big| ~var_{H}(\vx(k)) > t  ~\Big\} \\
 \mathcal{E}_{p}(t)=&\Big\{~k~\Big| ~ var_{H}(\vx(k)) \leq t ~\Big\}.
 \end{align*}
 and $ |\mathcal{T}| \leq Cs$ for some positive constant C. 
Since 
$
var(\bbE[\vx(k)|y])> \frac{C}{s^{\omega}},
$
we may choose  $t=\frac{a}{s^{\omega}}$ for sufficiently small positive constant $a$ such that for any $k \in \mathcal{T}$, $t<\frac{1}{2}var(\bbE[\vx(k)|y]$.
 According to Lemma \ref{lem:new_statistic:deviation} and the Bonferroni's inequality, we have
\begin{equation*}
\begin{aligned}
\bbP\left(\mathcal{T}^{c} \subset \mathcal{E}_{p}(t)\right)
&\geq 1-\sum_{k\in \mathcal{T}^{c}} \bbP\left(~var_{H}(\vx(k))> t\right) \\
&\geq 1-C_{1}\exp\left(-C_{2}\frac{ct}{H}+C_{3}\log(H)+\log(p-s)\right).
\end{aligned}
\end{equation*}
and  
\begin{equation*} 
\begin{aligned}
\bbP(\mathcal{T} \subset I_{p}(t))
&\geq\bbP\left(\bigcap_{k\in \mathcal{T} } \Big\{~var_{H}(\vx(k))\geq \frac{1}{2}var(\vm(k,y))~\Big\}\right) \\
&\geq 1-\sum_{k\in \mathcal{T}} \bbP\left(var_{H}\left(\vx(k)\right)<\frac{1}{2}var\left(\vm(k,y)\right)\right)\\
&\geq 1-C_1\exp\left(-C_2\frac{c~var(\vm(k,y))}{H}+C_{3}\log(H)+\log\left(Cs \right)\right),
\end{aligned}
\end{equation*}
$i.e.,$ we have  \eqref{ex_prob} and \eqref{in_prob} hold.
\epf

\subsection{Proof of Theorem \ref{thm:conditional_variance:high_dimensional}}
By choosing $H, c$ and $t=\frac{a}{s^{\omega}}$ properly, from Theorem \ref{thm:screen:consistent}, we have
\[
P(\widehat{\mathcal{T}}=\mathcal{T}) \geq 1-C_{1}\exp\left( -C_{2}\frac{n}{H^{2}s^{\omega}}+C_{3}\log(H)+\log(p-s) \right)
\]
for some positive constants $C_{i},i=1,2,3$. When $\widehat{\mathcal{T}}=\mathcal{T}$, we have  $|\widehat{\mathcal{T}}|=O(s)$. For the $n$ samples $(Y_{i},X^{\widehat{\mathcal{T}}}_{i})$, apply Theorem \ref{thm:high:dim}, we have
\begin{align*}
\|e(\widehat{\bLambda}_{H}^{\mathcal{T},\mathcal{T}})-\bLambda_p \|_{2}\leq
\|\widehat{\bLambda}_{H}^{\mathcal{T},\mathcal{T}}-\bLambda_{p}^{\mathcal{T},\mathcal{T}} \|_{2} \leq  O_{P}(
\frac{1}{H^{\vartheta}}+\frac{H^{2}s}{n}+\sqrt{\frac{H^{2}s}{n}}).
\end{align*} 
In particular, with probability converging to one, we have 
\[
\|e(\widehat{\bLambda}_{H}^{\mathcal{T},\mathcal{T}})-\bLambda_p \|_{2} \rightarrow 0 \mbox{ as } n \rightarrow \infty.
\]

\subsection{Proof of Theorem \ref{thm:operator_norem:high_dimensional}.}
The proof is almost identical to the proof of Theorem \ref{thm:space:consistency}, except that we additionally need to use Theorem 1 in \cite{bickel2008covariance}.


\section{Appendix C}\label{sec:assist:lemma}
\subsection{Assisting Lemmas}
\begin{definition}\label{def:treated.i.i.d.}
A set of random variables $x_{1},...,x_{n}$ can be treated as i.i.d random samples from a random variable $x$, if for any n variates symmetric function $f(w_{1},...,w_{n})$ , $f(x_{1},...,x_{n})$ is identically distributed as $f(z_{1},...,z_{n})$ where $z_{1},...,z_{n}$ are i.i.d random samples from $x$.
\end{definition}

\begin{lem}\label{rmk:key_rmk}
Let $(x_{i},y_{i})$ be n i.i.d random samples from a joint distribution $(x,y)$. Sort these samples according to the order statistics  of $y_i$'s and denote the sorted samples by $(x_{(1)}, y_{(1)}), (x_{(2)}, y_{(2)}), ..., (x_{(n)}, y_{(n)})$. Then for any a, b $(1\leq a \leq b+1\leq n)$, $x_{(a+1)},...,x_{(b)}$ can be treated   as $b-a$ i.i.d samples from $x\Big|\left(y\in [y_{(a)},y_{(b+1)}]\right)$ .

\proof In fact, we only need to prove that $y_{(a+1)}, \cdots,y_{(b)}$ can be treated as $b-a$ i.i.d. samples of $y\Big|(y\in [y_{(a)},y_{(b+1)}])$. The latter only needs to be proved for uniform distribution which  can be verified directly.
\end{lem}

\begin{cor}\label{cor:i.i.d}
In the slicing inverse regression contexts, recall that $S_{h}$ denotes the h-th interval $(y_{h-1,c},y_{h,c}]$ for $2\leq h \leq H-1$ and $S_{1}=(-\infty, y_{1,c}]$, $S_{H}=(y_{H-1,c},\infty)$. We have that
$x_{h,i}, i=1,\cdots,c-1$ can be treated as $c-1$  random samples of $x\Big|(y \in S_{h})$ for $h=1,...,H-1$ and $x_{H,1},...,x_{H}$ can be treated as $c$ random samples of $x\Big|(y\in S_{H})$. 
\end{cor}

\begin{lem}\label{lem:fourth_moment:lemma} 
Suppose that $(x,y)$ are defined over $\sigma$-finite space $\mathcal{X}\times\mathcal{Y}$ and g is a non-negative function such that $\bbE[g(x)]$ exists.  For any fixed positive constants $C_{1}<1<C_{2}$, there exists a constant C which only depends on $C_{1}, C_{2}$ such that for any partition 
$\mathbb{R}=\bigcup_{h=1}^{H}S_{h} $ where $S_{h}$ are intervals  satisfying
\begin{align}\label{partition:condition}
  \frac{C_{1}}{H} \leq \bbP(y \in S_{h})  \leq \frac{C_{2}}{H}, \forall h,
\end{align}
we have
\[
\sup_h\bbE(g(x)\Big|{y\in S'_{h}}) \leq  CH\bbE[g(x)].
\]
\proof
According to Fubini's Theorem, for any $h$,
\begin{align*}\label{eqn:Fubini:fourth:moments}
\bbE[g(x)]
=& \sum_{k}\bbP(y\in S_{k})\int_{\mathcal{X}} g(x)p(x|y\in S_{k})dx\\
\geq & \bbP(y\in S_{h})\int_{\mathcal{X}} g(x)p(x|y\in S_{h})dx.
\end{align*}
Due to the condition \eqref{partition:condition},  there exists a positive constant C such that
\[
\int_{\mathcal{X}} g(x)p(x|y\in S_{h})dx \leq CH\bbE[g(x)].
\]

\epf
\end{lem}

\begin{cor}\label{cor:restricted:covariance}
Let $\vx$ be a multivariate random variable with covariance matrix $\bSigma$. For any partition satisfying (\ref{partition:condition}), there exists a constant C such that
\[
var(\bbeta^{\tau}\vx|_{y\in S'_{h}})\leq CHvar(\bbeta^{\tau}\vx),\textrm{for any unit vector $\bbeta$},
\]
and
\[
\lambda_{max}\left( var\left(\vx\Big|{y\in S'_{h}}\right) \right) \leq CH\lambda_{max}\left( var\left(\vx\right) \right).
\]
\end{cor}

\begin{cor}\label{cor:restricted:sub-Gaussian}
Let $x$ be a sub-Gaussian random variable which is upper-exponentially bounded by K. Then for any partition satisfying (\ref{partition:condition}), there exists a constant C such that
\[
\bbE[\exp\left( \frac{x^{2}}{K^{2}}\right)\Big|{y\in S'_{h}}]\leq CH\bbE[\exp\left( \frac{x^{2}}{K^{2}}\right)].
\]
\end{cor}


Recall the definition of the random intervals $S_h, h=1,2,\cdots, H$ and random variable $\delta_h=\delta_h(\omega)=\int_{y\in S_h(\omega)}f(y)dy$.

\begin{lem} \label{lem:elementary:deviation}
Define the event $E(\epsilon)=\Big\{~\omega ~\Big|~ |\delta_{h}-\frac{1}{H}| > \epsilon, \forall h~\Big\}$. There exists a positive constant $C$ such that, for any  $\epsilon>\frac{4}{Hc-1} $ we have
\begin{equation}
P(E(\epsilon))\leq CH^{2}\sqrt{Hc+1}\exp(-(Hc+1)\frac{\epsilon^{2}}{32})
\end{equation}
for sufficient large H and  c. 
\proof The proof is deferred to the end of this paper.
\end{lem}

\subsection{Some Results from Random Matrices Theory.} We collect some direct corollaries of the non-asymptotic random matrices theory (e.g., \cite{rudelson2013hanson}).
\begin{lem} \label{thm:isometry_bound:sub_Gaussian}
Let $\vM$ be any $p\times n$ matrix $(n>p)$ whose columns $\vM_{i}$ are independent sub-Gaussian random vectors  in $\mathbb{R}^{p}$ with second moment $\vI_{p}$ and $\lambda_{sing,min}^{+}(\vM)$, $\lambda_{sing,max}(\vM)$ be the minimal non-zero and maximal singular value of $\vM$. Then for every t, with  probability at least $1-2\exp(-C't^{2})$, we have : 
\[
\sqrt{n}-C\sqrt{p}-t\leq \lambda^{+}_{sing,min}(\vM)\leq \lambda_{sing,max}(\vM)\leq \sqrt{n}+C\sqrt{p}+t.
\]
\end{lem}

\begin{lem}\label{lem:folklore}
Let $\vx_{1},\cdots,\vx_{n}$ be n i.i.d. samples from a p-dimensional sub-Gaussian random variable with covariance matrix $\bSigma$ and $\rho=\frac{p}{n}$. If there exists positive constants $C_1$ and $C_2$ such that 
\[
C_1\leq \lambda_{\min}(\bSigma_{\vx})\leq \lambda_{\max}(\bSigma_{\vx})\leq C_2.
\]
Let $\widehat{\bSigma}_{\vx} = \frac{1}{n}\sum_i \vx_i\vx_i^{\tau}$. Then 
\[
\|\widehat{\bSigma}_{\vx}-\bSigma_{\vx}\|_{2}\rightarrow 0 \mbox { if } \rho=0 \mbox{ when } n \rightarrow \infty.
\]
It is also easy to see that, given 
the boundedness condition on $\bSigma_{X}$ , $\|\widehat{\bSigma}^{-1}_{X}-\bSigma^{-1}_{X}\|_{2}\rightarrow 0$ if $\rho=\frac{p}{n}\rightarrow 0$ when $n\rightarrow \infty$.

\proof
Let $\vx_{i}=\bSigma_{\vx}^{1/2}\vm_{i}$ where $\vm_{i}$ is sub-Gaussian random variable with covariance  matrix $\vI_{p}$ and $\vM=(\vm_{1},...,\vm_{n})$.
From Lemma \ref{thm:isometry_bound:sub_Gaussian}, we have
\begin{align*}
&\|\frac{1}{n}\vM\vM^{\tau}-\vI_{p}\|_{2}\rightarrow 0\\
\mbox{ \noindent and \quad\quad\quad \quad \quad\quad \quad} &\\
\|\widehat{\bSigma}^{-1}-\bSigma^{-1}\|_{2}&=\|\bSigma^{-1/2}\|_{2}\|\frac{1}{n}\vM\vM^{\tau}-\vI_{p}\|_{2}\|\bSigma^{-1/2}\|_{2}\rightarrow 0,
\end{align*}
with probability converges to 1 as $n \rightarrow \infty$ because
\begin{align*}
\lambda_{max}\left(\frac{1}{n}\vM\vM^{\tau}\right)&\leq \left(1+\frac{(C+1)\sqrt{p}}{\sqrt{n}}\right)^{2}\\
\mbox{and}\\
\lambda_{min}\left(\frac{1}{n}\vM\vM^{\tau}\right)&\geq \left(1-\frac{(C+1)\sqrt{p}}{\sqrt{n}}\right)^{2}
\end{align*}
with probability at least $1-2\exp(-C'p)$.
\epf
\end{lem}

The following lemma is well known in the non-asymptotic random matrix theory (\cite{Vershynin:2010} Proposition 5.34 ) which is slightly different from the Lemma \ref{thm:isometry_bound:sub_Gaussian}.
\begin{lem} \label{random:nonasymptotic}
Let $\vE_{p\times H}$ be a $p\times H$ matrix, whose entries are independent standard normal random variables. Then for every $ t\geq 0$, with probability at least $1-2\exp(-t^{2}/2)$, we have :
\[ \lambda^{+}_{sing,min}(\vE_{p\times H}) \geq \sqrt{p}-\sqrt{H}-t,\]
and
\[ \lambda_{sing,max}(\vE_{p\times H}) \leq \sqrt{p}+\sqrt{H}+t.\]
\end{lem} 

\begin{cor}\label{cor:final:done}
We have
\[
\frac{1}{2}\left(\sqrt{p}-\sqrt{H}\right) \leq \lambda^{-}_{sing,min}\left(\vE_{p\times H}\right) \leq \lambda_{sing,max}\left(\vE_{p\times H}\right)\leq \frac{3}{2}\left(\sqrt{p}+\sqrt{H}\right).
\]
with probability converging to one, as $n \rightarrow \infty$.
\proof

Choosing $t=\sqrt{p}/2$, according to Lemma \ref{random:nonasymptotic}, we have:
\[
\bbP\left(\frac{\lambda_{max}(E_{H})}{\sqrt{p}+\sqrt{H}} \leq \frac{3}{2}\right) \geq \bbP\left(\frac{\lambda_{max}(E_{H})}{\sqrt{p}+\sqrt{H}} \leq 1+\frac{\sqrt{p}}{2\sqrt{p}+2\sqrt{H}}\right) 
\]
and 
\[
\bbP\left(\frac{\lambda^{+}_{min}(E_{H})}{\sqrt{p}-\sqrt{H}} \geq \frac{1}{2} \right) \geq  \bbP\left(\frac{\lambda_{max}(E_{H})}{\sqrt{p}-\sqrt{H}} \geq 1-\frac{\sqrt{p}}{2\sqrt{p}-2\sqrt{H}}\right) 
\]
with probability at least $1-2\exp(-p/8)$.
i.e., With probability converging to one, we have
\[
\frac{1}{2}(\sqrt{p}-\sqrt{H}) \leq \lambda^{-}_{min}(E_{p\times H}) \leq \lambda_{max}(E_{p\times H})\leq \frac{3}{2}(\sqrt{p}+\sqrt{H}).
\]
 \epf
\end{cor}

\subsection{Basic Properties of sub-Gaussian random variables.}

We rephrased several equivalent definitions of the sub-Gaussian distribution here (See e.g., \cite{Vershynin:2010} ):
\begin{definition} \label{def:subGaussian:definition}
Let x be a random variable. Then the following properties are equivalent with parameters $K_{i}$'s differing from each other by at most an absolute constant factor,
\begin{itemize}
\item [1.] Tails: $\bbP(|x|>t) \leq \exp(1-t^{2}/K^{2}_{1})$ for all $t\geq 0$.
\item [2.] Moments: $(\bbE[|x|^{p}])^{1/p} \leq K_{2}\sqrt{p}$ for all $p\geq 1$.
\item [3.] Super-exponential moment: $\bbE\exp(x^{2}/K^{2}_{3})\leq e$.
\end{itemize}
Moreover, if $\bbE[x]=0$, then the properties $1-3$ are also equivalent to the following one:
\begin{itemize}
\item [4.] Moment generating function: $\bbE[\exp(tx)] \leq \exp(t^{2}K^{2}_{4})$.
\end{itemize}
\end{definition}

\begin{definition}\label{def:subgaussion:upper:norm}
For a sub-Gaussian random variable $x$ with the constants $K_{i},i=1,2,3,4$ given in Definition \ref{def:subGaussian:definition},
we will call a constant $K$ an upper-exponential bound  of $x$ or $x$ is upper-exponentially bounded by  $K$ if $K> \max_{i}\{K_{1}, K_{2},K_{3},K_{4}\}$.
\end{definition}  
We summarize some properties regarding the sub-Gaussian distributions into the following lemmas.

\begin{lem} \label{lem:subGaussian:tail}  Let $\delta_{1},...,\delta_{n}$ be $n$ $($not necessarily independent or with mean zero$)$ sub-Gaussian random variables upper-exponentially bounded by K.
\begin{itemize}
\item[{\bf i)}]  $\frac{1}{n}\sum_{i=1}^{n}\delta_{i}$ is sub-Gaussian and upper-exponentially bounded by K.

\vspace{3mm}
\item[{\bf ii)}] $\delta_{1}-\bbE[\delta_{1}]$ is sub-Gaussian upper-exponentially bounded by 2K.
\item[{\bf iii)}] If they are  independent and with mean zero , then $\frac{1}{\sqrt{n}}\sum_{i}\delta_{i}$ is sub-Gaussian and upper-exponentially bounded by K. 

\vspace{3mm}
\item[{\bf iv)}] If they are  i.i.d., then we have the concentration inequality:
\[
\bbP\left(\big|\frac{\sum_{i=1}^{n}\delta_{i}}{n}-\bbE[x]\big|>t\right) \leq 2\exp\left( \frac{-nt^{2}}{2K^{2}e+2tK} \right) .
\]

\end{itemize}

\proof

{\noindent\bf i) } follows from the linear property of expectation and the the convexity of exponential function. i.e.,
\begin{align*}
\bbE[\exp(\frac{1}{nK^{2}}\sum_{i}\delta^{2}_{i})]\leq \bbE[\frac{1}{n}\sum_{h}\exp(\frac{\delta^{2}_{i}}{K^{2}})] \leq \max_{i}\bbE[\exp(\frac{\delta^{2}_{i}}{K^{2}})] \leq e.
\end{align*}

{\noindent\bf ii) } From Definition \ref{def:subGaussian:definition}, we know that $|\bbE[\delta_{i}] | \leq K$ which gives us the desired upper-exponential bound of $\delta_{i}-\bbE[\delta_{i}]$.

{\noindent\bf iii)} is trivial as $\delta_{1},\cdots,\delta_{c}$ are independent and with mean zero.

{\noindent \bf iv)} Since $\delta_{1}$ is sub-Gaussian upper-exponentially bounded by K, we have:
\begin{align}
\nonumber  \bbE[|\delta_{1}|^{p}] &= \int_{0}^{\infty} pt^{p-1}\bbP(|\delta_{1}|>t)dt \leq \int_{0}^{\infty}pt^{p-1}\exp\left(1-\frac{t^{2}}{K^{2}} \right)dt\\
\nonumber &  = \frac{ep}{2}\Gamma(\frac{p}{2})K^{p} \quad &\mbox{ for any  $p\geq 1$}     \\ 
\nonumber & \leq p!K^{p-2} \frac{(K^{2}e)}{2}  \quad &\mbox{ for any $p\geq 2$}
\end{align}

Recall the well known Bernstein inequality.
\begin{lem} {\bf ( Bernstein Inequality ).}
If there exists positive constants V and b such that for any integers $m \geq 2$, 
\[
E[|\delta_{1}|^{m}] \leq m!b^{m-2}V/2
\]

then
\begin{equation}\label{eqn:bernstein:inequality}
\bbP\left(\big|\frac{\sum_{i=1}^{n}\delta_{i}}{n}-\bbE[x]\big|>t\right) \leq 2\exp\left( \frac{-nt^{2}}{2V+2tb} \right) .
\end{equation}
\end{lem}
By chooing $V=K^{2}e$ and $b=K$, we get the desired concentration inequality. \epf
\end{lem}

\begin{lem}\label{cor:key_deviation:lemma} 
Suppose that $(x,y)$ are defined over $\sigma$-finite space $\mathcal{X}\times\mathcal{Y}$ and x is sub-Gaussian with mean 0 and upper exponentially bounded by $K$, let $m(y)=\bbE[x|y]$, $\epsilon(y)=x-m(y)$, then we have
\begin{itemize}
\vspace*{3mm}
\item[{\bf i)}]  $m(y)$ and $\epsilon	(y)$ are sub-Gaussian and upper-exponentially bounded by $K$ and $2K$ respectively.

\vspace*{3mm}
\item[{\bf ii)}] Let $\mathcal{Z}$ consists of points y such that $x|_{y}$ is not sub-Gaussian, i.e.,
\[
\mathcal{Z}\triangleq\Big\{y~| ~\exists t \in (0,t_{0}] \mbox{ such that } \int_{\mathcal{X}} \exp(tx^{2})p(x|y)p(y)dx =\infty \Big\},
\]
then   $\bbP(y\in \mathcal{Z})=0$. 

\vspace*{3mm}
\item[{\bf iii)}] For any fixed positive constants $C_{1}<1<C_{2}$ and any partition 
$\mathbb{R}=\bigcup_{h=1}^{H}S_{h} $ where $S_{h}$ are intervals
 satisfying
 \[
\frac{C_{1}}{H} \leq \bbP(y \in S_{h})  \leq \frac{C_{2}}{H}, \forall h,
 \]
  there exists a constant C such that
 \[
\sup_h\bbP(x|_{y\in S_{h}}>t) \leq  CH\exp\left(1-\frac{t^{2}}{K^{2}} \right) 
 .\]
As a direct corollary, we know that there exists a positive constant C such that
\[
\bbE\left[\exp\left( \frac{(x|_{y\in S_{h}})^{2}}{2K^{2}}\right)\right] \leq CH,
\]
and 
\[
\bbE\left[\Big|(x|_{y\in S_{h}})\Big|^{m}\right] \leq CHmK^{m}\Gamma(\frac{m}{2})/2.
\]

\vspace*{3mm}
 \item[{\bf vi)}] Suppose that $x|_{y \in S_{h}}$ is defined as in ${\bf iii)} $. Let $x_{i}$, i=1,...,c be c samples from $x|_{y\in S_{h}}$, $\bar{x}_{h}=\frac{1}{c} \sum_{i}x_{i}$ and $\mu_{h}=\bbE[x|_{y\in S_{h}}]$, we have
 \[
\bbP[|\bar{x}_{h}-\mu_{h}| > t]\leq 2 \exp\left( \frac{-ct^{2}}{2CHK^{2}+2tK}\right).
 \]

 \end{itemize}
\proof 
 {\bf i)}  By Jensen's inequality, we have
\begin{align*}
\bbE[\exp(t\bbE[x|y])] \leq \bbE[\bbE[\exp(tx)|y]]=\bbE[\exp(tx)]\leq \exp(t^{2}K_{1}^{2}).
\end{align*}
i.e., $m(y)$ is sub-Gaussian and upper-exponentially bounded by $K_{1}$. Since $x$ , $m(y)$ is sub-Gaussian and upper-exponentially bounded by $K_{1}$, we know that $\epsilon=x-m(y)$ is sub-Gaussian and upper-exponentially bounded by $2K_{1}$.

\vspace{3mm}
 \item [{\bf ii)}] Let $p(x,y)$ be the joint density function of $(x,y)$ and $p(x)$, $p(y)$ be the marginal distribution of x, y. Since x is sub-Gaussian, we know there exists $t_{0}>0$ such that
 \[
\int_{\mathcal{X}}\exp(tx^{2})\int_{\mathcal{Y}} p(x|y)p(y)dy dx \leq e	\mbox{ for } 0\leq t\leq t_{0}.
 \]
By Fubini Theorem, we know
\begin{equation}\label{eqn:Fubini}
\int_{\mathcal{Y}}p(y)\int_{\mathcal{X}} \exp(tx^{2})p(x|y)dxdy \leq e \mbox{ for } 0\leq t\leq t_{0}.
\end{equation}

Recall that we have $\mathcal{Z}\triangleq\{y| \exists t \in (0,t_{0}] \mbox{ such that } \int_{\mathcal{X}} \exp(tx^{2})p(x|y)p(y)dx =\infty \}$, from equation \eqref{eqn:Fubini}, we know $\bbP(y\in \mathcal{Z})$=0. In particular, we know that for any $y \not \in \mathcal{Z}$, $x|_{y}$ is sub-Gaussian. However, the norm (e.g., sub-exponential norm)of $x|_{y}$ might be varying along with y and , as a function of y , it might be not bounded. 

\vspace*{3mm}
\item[{\bf iii)}] From Lemma \ref{lem:fourth_moment:lemma},
we know that  there exists a positive constant C such that
\[
\int_{\mathcal{X}} \exp(tx^{2})p(x|y\in S_{h})dx \leq CHe.
\]
 For simplicity if notation, we will denote $x|_{y\in S_{h}}$  by $z$ through out this lemma. 
So $\mbox{ for } 0\leq t\leq t_{0}=\frac{1}{K}$, we have
\[
\bbP\left( z>a\right) \leq \frac{\bbE[\exp\left( tz)^{2}\right)]}{\exp(t^{2}a^{2})} \leq CHe\exp\left(-t^{2}a^{2} \right).
\]
From the above tail bounds, we have that for any integer $m>0$
\begin{align*}
\bbE[\big|z\big|^{m}]=&\int_{0}^{\infty}\bbP(|z|>t)(m)t^{m-1}dt\leq CHm\int_{0}^{\infty}exp(-\frac{t^{2}}{K^{2}})t^{m-1}dt\\
\leq & CHmK^{m}(m/2)/2.
\end{align*}
We then have
\begin{align*}
\bbE\left[ \exp\left( tz^{2} \right)\right] &\leq \sum_{m=0}^{\infty}\frac{\bbE\left[t^{m} z^{2m}\right]}{m!} \leq \sum_{m=0}^{\infty}\frac{\bbE\left[t^{m}z^{2m}\right]}{m!}\\
&\leq \sum_{m=0}^{\infty}\frac{t^{m}CHmK^{2m}\Gamma(m)}{m!}=CH\sum_{m=0}^{\infty}t^{m}K^{2m}
\end{align*}
From which we know if $0\leq t<\frac{1}{2}K^{-2}$, the R.H.S is bounded by $CH$ for a positive constant C.

\vspace*{3mm}
\item[{\bf vi)}] From the previous proof, we know that for any integer $m\geq 2$
\[
\bbE[|z|^{m}] \leq CHm!K^{m}=m!K^{m-2} (2CHK^{2})/2.
\]
By the Bernstein inequality \eqref{eqn:bernstein:inequality}, we have:
\[
\bbP\left(\Big|\frac{\sum_{i=1}^{c}z_{i}}{c}-\bbE[z]\Big|>t\right) \leq 2 \exp\left( \frac{-ct^{2}}{2CHK^{2}+2tK}\right).
\]

\epf
\end{lem}

\begin{lem} \label{lem:variance_deviation}
Let $z_{i}, i=1,\cdots,n$ be i.i.d. samples of a sub-Gaussian distribution exponentially upper bounded by K,  
then there exist positive constants $C_{1}, C_{2}$ such that, if $\sqrt{n}\epsilon\rightarrow \infty$, we have
\[
\mathbb{P}(|\frac{1}{n}\sum_{i}(z_{i}-\bar{z})^{2}-var(z)|>\epsilon)\leq C_{1}\exp(-C_{2} \frac{\epsilon \sqrt{n}}{K^{2}}),
\]
where $\bar{z}=\frac{1}{n}\sum_{i}z_{i}$.
\proof 
Recall the following Hanson-Wright inequality in \cite{rudelson2013hanson}

\begin{lem}\label{lem:Hanson_Wright}
Let $\vv=(\vx(1),\cdots,\vx(n))$ be a sub-Gaussian random vector with independent components $\vx(\bbeta)$ such that $\bbE[\vx(\bbeta)]=0$ and $\|\vx(\bbeta)\|_{\psi_{2}}\leq K$. Let $\vA$ be an $n\times n$ matrix. Then there exists a positive constant C such that for any $t>0$, 
\[
\mathbb{P}\{|\vx^{\tau}\vA\vx-\bbE[\vx^{\tau}A\vx]|>t\}\leq 2 \exp\left(-C ~\min\left(\frac{t^{2}}{K^{4}\|\vA\|^{2}_{HS}},\frac{t}{K^{2}\|\vA\|_{HS}}\right) \right).
\]
Here the $\psi_{2}$ norm of a random variable z is defined as 
$\|z\|_{\psi_{2}}\triangleq\sup_{p}p^{-1/2}(\bbE|z|^{p})^{1/p}$ and 
the HS norm of a matrix $\vA$ is defined as $\|\vA\|_{HS}=(\sum_{i,j}|a_{i,j}|^{2})^{1/2}$.
\end{lem} 

Since
\begin{align*}
&\bbP\left(|\frac{1}{n}\sum_{i}(z_{i}-\bar{z})^{2}-var(z)|> 2\epsilon\right) \\
&=\bbP\left(|\frac{1}{n}\sum_{i}(z_{i}-\bbE[z])^{2}-(\bbE[z]-\bar{z})^{2}-\bbE[(z-E[z])^{2}]|> 2\epsilon\right)\\
&\leq 
\bbP\left(| \frac{1}{n}\sum_{i}(z_{i}-\bbE[z])^{2}-\bbE[(z-\bbE[z])^{2}] |>\epsilon\right)+
\bbP\left((\bbE[z]-\bar{z})^{2}>\epsilon\right),
\end{align*}
and $z_{i}-\bbE[z]$ are sub-Gaussian with mean 0,  from Lemma \ref{lem:Hanson_Wright} by choosing $\vA=\frac{1}{n}\vI_{p}$ and $\vz^{\tau}=(z_{1}-\bbE[z],..,z_{p}-\bbE[z])$, we have
\begin{equation}\label{eqn:temp:no_names}
\mathbb{P}\left(| \frac{1}{n}\sum_{i}(z_{i}-\bbE[z])^{2}-\bbE[(z-\bbE[z])^{2}] |>\epsilon\right) \leq 2\exp\left(-C \frac{\sqrt{n}\epsilon}{K^{2}}\right),
\end{equation}
since $\sqrt{n}\epsilon \rightarrow \infty$. 

The following follows from the usual deviation argument:
\begin{align*}
\mathbb{P}\left((\bbE[z]-\bar{z})^{2}>\epsilon\right)
&\leq C_{1}\exp\left(-C_{2}n\epsilon\right).
\end{align*}

Combining with the estimate \eqref{eqn:temp:no_names},  we know there exists positive constants $C_{1}$ and $C_{2}$ such that
\[
\mathbb{P}\left(|\frac{1}{n}\sum_{i}(z_{i}-\bar{z})^{2}-var(z)|> \epsilon\right) \leq C_{1}\exp\left(-C_{2}\frac{\sqrt{n}\epsilon}{K^{2}}\right),
\]
for sufficiently large n since $\sqrt{n}\epsilon\rightarrow \infty$. \epf

\end{lem}

\subsection{ Proof of Lemma \ref{lem:elementary:deviation}. } We only need to prove this lemma for n i.i.d. sample $y_i$'s from a uniform distribution over $[0,1]$.  We slightly change the notation of order statistics $y_{(i)}$ to $y_{(i,n)}$ so that we can keep track of the sample size. Since y is uniform distribution on $[0,1]$, it is well known that $y_{(i,n)} \sim Beta(i,n-i+1)$ with expectation $\frac{i}{n+1}$ and mode $\frac{i-1}{n-1}$.  
Lemma \ref{lem:elementary:deviation} is a direct corollary of the following lemma. 

\begin{lem} \label{lem:unirom:deviation:right} Suppose there are $n=Hc$ i.i.d. samples from uniform distribution over $[0,1]$, when H,c are sufficiently large, we have the following large deviation inequalities of $y_{(kc,Hc)}$, $k=1,\cdots,(H-1)$.
\begin{itemize}
\item[i)] There exists a positive constant C, such that for any $\epsilon > \frac{1}{Hc-1}$, we have
\[
\bbP\left(y_{(kc,Hc)}>\frac{k}{H}+\epsilon\right) \leq CH\sqrt{Hc+1}\exp\left(-(Hc+1)\frac{\epsilon^{2}}{2}\right);
\] 
\item[ii)] 
 There exists a positive constant C, such that for any $\epsilon > \frac{2}{Hc-1} $, we have
\[
\bbP\left(y_{(kc,Hc)}<\frac{k}{H}-\epsilon\right)\leq CH\sqrt{Hc+1}\exp\left(-(Hc+1)\frac{\epsilon^{2}}{8}\right);
\]
\item[iii)] Let $\delta(k,H,c)=|y_{((k-1)c,Hc)}-y_{(kc,Hc)}|$ , for $2\leq k \leq H-1$, $\delta(1,H,c)=|y_{(c,Hc)}|$ and $\delta(H,H,c)=|1-y_{((H-1)c,Hc)}|$.  There exists a positive constant C, such that for any $\epsilon>\frac{4}{Hc-1} $  , we have
for any $1\leq k\leq H$:
\[
\bbP\left(|\delta(k,H,c)-\frac{1}{H}|>\epsilon\right) \leq CH\sqrt{Hc+1}\exp\left(-(Hc+1)\frac{\epsilon^{2}}{32}\right).
\]
\end{itemize}
\end{lem}
\vspace*{3mm}
We will prove Lemma \ref{lem:unirom:deviation:right} later. Assuming it, we have
\begin{align*}
\bbP\left( E(\epsilon) \right) \leq \sum_{k=1}^{H}\bbP\Big(|\delta(k,H,c)-\frac{1}{H}|&>\epsilon\Big)\\
\leq CH^{2}\sqrt{Hc+1}&\exp\left(-(Hc+1)\frac{\epsilon^{2}}{32}\right).
\end{align*}
  \epf

\subsubsection{ Proof of Lemma \ref{lem:unirom:deviation:right}.}  {\it The first part: } For any $1\leq k\leq H-1$, we note that
\begin{align*}
&\bbP\left(y_{(kc,Hc)}>\frac{k}{H}+\epsilon\right) \leq \bbP\left(y_{(kc,Hc)}>\frac{kc}{Hc+1}+\epsilon\right)\\
&=\frac{1}{B(kc,Hc-kc+1)}\int^{1}_{x>\frac{kc}{Hc+1}+\epsilon}x^{kc-1}(1-x)^{Hc-kc}dx.
\end{align*}
When $\epsilon>\frac{1}{Hc-1}$, we know
the mode $x_{M}=\frac{kc-1}{Hc-1}<x_{D}\triangleq\frac{kc}{Hc+1}+\epsilon$ , so we have

\begin{equation}\label{eqn:inline:lemma:1}
\begin{aligned}
\bbP\left(y_{(kc,Hc)}>\frac{k}{H}+\epsilon\right)&\leq \frac{(x_{D})^{kc-1}(1-x_{D})^{Hc-kc+1}}{B(kc,Hc-kc+1)}\\
&\leq H\frac{(x_{D})^{kc}(1-x_{D})^{Hc-kc+1}}{B(kc,Hc-kc+1)}.
\end{aligned}
\end{equation}
The last inequality due to $Hx_{D}\geq 1$.  If $\epsilon+\frac{k}{H}\geq 1$, then $\bbP(y_{(kc,Hc)}>\frac{k}{H}+\epsilon) =0$ and Lemma \ref{lem:unirom:deviation:right}  holds automatically. 

We may assume that $\epsilon+\frac{k}{H}< 1$ below. Let us denote the right hand side of  \eqref{eqn:inline:lemma:1} by A, then
\begin{align*}
\log(A)=&\log(H)+kc\log(E+\epsilon)+(Hc-kc+1)\log(1-E-\epsilon)\\
&+\log(Hc+1)!-\log(kc)!-\log(Hc-kc+1)! \\
&-\log(Hc+1)+\log(kc)+\log(Hc-kc+1),
\end{align*}
where $E=\frac{kc}{Hc+1}$. According to the Stirling formula:
\[
\log(n!)=n\log(n)-n+\frac{1}{2}\log(2\pi n)+O(\frac{1}{n}),
\]
when $m$ is sufficiently large we have:
\begin{align}
\log(A) =&\log(H)+ (Hc+1)(E\log(1+\frac{\epsilon}{E})+(1-E)\log(1-\frac{\epsilon}{1-E})) \nonumber \\
&-\frac{1}{2}(\log(Hc+1)-\log(kc)-\log(Hc-kc+1))\nonumber \\
&-\frac{1}{2}\log(2\pi)+O(\frac{1}{kc})+O(\frac{1}{Hc-kc+1})\nonumber \\
\leq&\log(H) -\frac{(Hc+1)\epsilon^{2}}{2(1-E)}-\frac{1}{2}\log(2\pi)  \label{eqn:inline:temp:1} \\ 
&-\frac{1}{2}(\log(Hc+1)-\log(kc)-\log(Hc-kc+1))+O(\frac{1}{c}) \nonumber \\
\leq &\log(H)-\frac{(Hc+1)\epsilon^{2}}{2(1-E)}
-\frac{1}{2}(\log(Hc+1)-\log(kc)-\log(Hc-kc+1)) \nonumber ,
\end{align}
where we use the fact that $\frac{kc}{Hc+1}\leq E+\epsilon<1 $ and the following elementary lemma, which can be proved by the Taylor expansion:
\begin{lem}\label{lem:elementary:log}
Suppose $a, b$ are positive numbers such that $a+b=1$, then for any $0<\epsilon<b$, we have:
\[
a\log(1+\frac{\epsilon}{a})+b\log(1-\frac{\epsilon}{b}) \leq -\frac{\epsilon^{2}}{2b}.
\]
\end{lem}
 Now we know that there exists a positive constant C such that for any $1\leq k\leq H-1$ and for any $\epsilon> \frac{1}{Hc-1} $, the following holds:
\begin{align*} \label{bound:right}
\bbP\Big(y_{(kc,Hc)}&>\frac{k}{H}+\epsilon\Big) \\
&\leq CH\sqrt{\frac{(kc)(Hc-kc+1)}{Hc+1}}\exp\left(-(Hc+1)\frac{\epsilon^{2}}{2(1-E)}\right)\\
&\leq CH\sqrt{Hc+1}\exp\left(-(Hc+1)\frac{\epsilon^{2}}{2(1-E)}\right)\\
&\leq C H\sqrt{Hc+1}\exp\left(-(Hc+1)\frac{\epsilon^{2}}{2}\right).
\end{align*}
The last inequality follows from $\frac{\epsilon^{2}}{1-E} \geq \epsilon^{2}$ since $\frac{1}{H+1}\leq E\leq \frac{H-1}{H}$. 

\vspace{3mm}
{\it The second part: } The proof of the second part is similar. For completeness, we sketch some  calculations below. For any $1\leq k\leq H-1$, when $\epsilon>\frac{2}{Hc-1}$, we have
\begin{align*}
&\bbP\left(y_{(kc,Hc)}<\frac{k}{H}-\epsilon\right) \leq \bbP\left(y_{(kc,Hc)}<\frac{kc}{Hc+1}-\epsilon/2\right)   \\
&=\frac{1}{B(kc,Hc-kc+1)}\int_{x<\frac{kc}{Hc+1}-\epsilon}x^{kc-1}(1-x)^{Hc-kc}dx.
\end{align*}
Since $\epsilon>\frac{1}{Hc-1}$, we know
the mode $x_{M}=\frac{kc-1}{Hc-1}>x_{D'}\triangleq \frac{kc}{Hc+1}-\epsilon/2$, so we have
\begin{align*}
\bbP\left(y_{(kc,Hc)}\leq \frac{k}{H}-\epsilon\right)\leq& \frac{(x_{D'})^{kc}(1-x_{D'})^{Hc-kc}}{B(kc,Hc-kc+1)}\\
\leq& H\frac{(x_{D'})^{kc}(1-x_{D'})^{Hc-kc+1}}{B(kc,Hc-kc+1)}.
\end{align*}
The last inequality due to $H(1-x_{D'})\geq 1$.

The rest is similar to the first part. We have that for any $1\leq k\leq H-1$ and for any $ \epsilon >\frac{2}{Hc-1}$,
\begin{align}
\bbP\left(y_{(kc,Hc)}<\frac{k}{H+1}-\epsilon\right) &\leq CH\sqrt{Hc+1}\exp\left(-(Hc+1)\frac{\epsilon^{2}}{8}\right).
\end{align}

{\it The third  part: \quad } The third part is a direct corollary of the first two parts. Note that for any $2\leq k\leq H-2$, for any $\epsilon>\frac{4}{Hc-1}$
\begin{align*}
\bbP\Big(|\delta(k,H,c)&-\frac{1}{H}|>\epsilon\Big) \\
&=
\bbP\left(\Big | y_{(k+1)c,Hc}-\frac{k+1}{H}-(y_{kc,Hc}-\frac{k}{H}) \Big| >\epsilon\right)\\
&\leq\bbP\left(|y_{(k+1)c,Hc}-\frac{k+1}{H}|>\frac{\epsilon}{2}\right)+\bbP\left(|y_{kc,Hc}-\frac{k}{H}|>\frac{\epsilon}{2}\right)\\
&\leq CH\sqrt{Hc+1}\exp\left(-(Hc+1)\frac{\epsilon^{2}}{32}\right).
\end{align*}
When $k=1$, we have
\begin{align*}
\bbP\Big(|\delta(1,H,c)&-\frac{1}{H}|>\epsilon\Big) =
\bbP\left(|y_{c,Hc}-\frac{1}{H}|>\epsilon\right)\\
&\leq\bbP\left( y_{c,Hc}-\frac{1}{H}>\epsilon\right)+\bbP\left(y_{c,Hc}-\frac{1}{H}<-\epsilon\right)\\
&\leq CH\sqrt{Hc+1}\exp\left(-(Hc+1)\frac{\epsilon^{2}}{8}\right)\\
&\leq CH\sqrt{Hc+1}\exp\left(-(Hc+1)\frac{\epsilon^{2}}{32}\right).
\end{align*}
When $k=H$, we have
\begin{align*}
\bbP\Big( |\delta(H-1,H,c)&-\frac{1}{H}|>\epsilon \Big)
=
\bbP\left(|y_{(H-1)c,Hc}-\frac{H-1}{H}|>\epsilon\right)\\
&\leq\bbP\left(y_{(k+1)c,Hc}-\frac{H-1}{H}>\epsilon\right)+\bbP\left(y_{(H-1)c,Hc}-\frac{H-1}{H}<-\epsilon\right)\\
&\leq CH\sqrt{Hc+1}\exp\left(-(Hc+1)\frac{\epsilon^{2}}{8}\right)\\
&\leq CH\sqrt{Hc+1}\exp\left(-(Hc+1)\frac{\epsilon^{2}}{32}\right).
\end{align*}
\epf

\end{appendices}

\bibliographystyle{plainnat}
\bibliography{sir}

\end{document}